\documentclass{llncs}

\usepackage{makeidx}  % allows for indexgeneration
\usepackage{amsfonts,amssymb,amsmath,amscd} %,amsmath,amsthm,amscd}
\usepackage{graphicx}
\usepackage{enumerate}

\newcommand{\ZZ}{\mathbb{Z}}
\newcommand{\NN}{\mathbb{N}}
\newcommand{\QQ}{\mathbb{Q}}
\newcommand{\RR}{\mathbb{R}}

\newcommand{\PPP}{\mathcal{P}}
\newcommand{\MMM}{\mathcal{M}}
\newcommand{\OOO}{\mathcal{O}}

\newcommand{\rar}{\rightarrow}

\newcommand{\mset}[2]{\left\{ #1 \; \middle| \; #2 \right\}}
\newcommand{\mmset}[2]{\{ #1 \; | \; #2 \}}
\newcommand{\choice}[1]{\left\{ \begin{array}{ll} #1 \end{array} \right.}
\newcommand{\floor}[1]{\left\lfloor #1 \right\rfloor}
\newcommand{\ceil}[1]{\left\lceil #1 \right\rceil}
\newcommand{\mmod}{ \; \operatorname{mod} \;}
\newcommand{\sprod}[2]{\left\langle #1 , #2 \right\rangle}
\newcommand{\mat}[1]{\begin{pmatrix}#1\end{pmatrix}}

\newcommand{\ehr}{\operatorname{ehr}}
\newcommand{\lcm}{\operatorname{lcm}}
\newcommand{\supp}{\operatorname{supp}}
\newcommand{\cone}{\operatorname{cone}}
\newcommand{\vcone}{\operatorname{vcone}}

\begin{document}

\pagestyle{headings}  % switches on printing of running heads
\mainmatter              % start of the contributions
\title{An Invitation to Ehrhart Theory: Polyhedral Geometry and its Applications in Enumerative Combinatorics}
\toctitle{An Invitation to Ehrhart Theory: Polyhedral Geometry and its Applications in Enumerative Combinatorics}
\titlerunning{Invitation to Ehrhart Theory}  

\author{Felix Breuer\thanks{Felix Breuer was supported by Austrian Science Fund (FWF) special research group \emph{Algorithmic and Enumerative Combinatorics} SFB F50-06.}}
\authorrunning{Felix Breuer} % abbreviated author list (for running head)
\institute{Research Institute for Symbolic Computation\\
Johannes Kepler University\\
Altenberger Str. 69,
4040 Linz, Austria\\
\emph{email:} \email{felix@felixbreuer.net},\\
\emph{web:} \texttt{http://www.felixbreuer.net}}

\maketitle              % typeset the title of the contribution

\begin{abstract}
In this expository article we give an introduction to Ehrhart theory, i.e., the theory of integer points in polyhedra, and take a tour through its applications in enumerative combinatorics. Topics include geometric modeling in combinatorics, Ehrhart's method for proving that a counting function is a polynomial, the connection between polyhedral cones, rational functions and quasisymmetric functions, methods for bounding coefficients, combinatorial reciprocity theorems, algorithms for counting integer points in polyhedra and computing rational function representations, as well as visualizations of the greatest common divisor and the Euclidean algorithm.

\end{abstract}

\section{Introduction}
\label{sec:intro}

Polyhedral geometry is a powerful tool for making the structure underlying many combinatorial problems visible -- often literally! In this expository article we give an introduction to Ehrhart theory and more generally the theory of integer points in polyhedra and take a tour through some of its many applications, especially in enumerative combinatorics.

In Section~\ref{sec:modeling}, we start with two classic examples of geometric modeling in combinatorics and then introduce Ehrhart's method for showing that a counting function is a (quasi-)polynomial in Section~\ref{sec:ehrhart}. We present combinatorial reciprocity theorems as a first application in Section~\ref{sec:reciprocity}, before we talk about cones as the basic building block of Ehrhart theory in Section~\ref{sec:cones}. The connection of cones to rational functions is the topic of Section~\ref{sec:brion}, followed by methods for proving bounds on the coefficients of Ehrhart polynomials in Section~\ref{sec:coefficients}. Section~\ref{sec:quasisymmetric} discusses a surprising connection to quasisymmetric functions. Section~\ref{sec:computing} is about algorithms for counting integer points in polyhedra and computing rational function representations, in particular Barvinok's theorem on short rational functions. Finally, Section~\ref{sec:euclid} closes with a playful look at the connection between the Euclidean algorithm and the geometry of $\ZZ^2$.

\section{Geometric Modeling in Combinatorics}
\label{sec:modeling}

Many objects in combinatorics can be conveniently modeled as integer vectors that satisfy a set of linear equations and inequalities. In applied mathematics, this paradigm has proven tremendously successful: the combinatorial optimization industry rests to a large part on mixed integer programming. However, also in pure mathematics this approach can help to prove theorems. We illustrate this approach of constructing geometric models of combinatorial objects and problems on two of the most classic counting functions in all of combinatorics: The chromatic polynomial of a graph and the restricted partition function.

\begin{figure}[t]
\includegraphics[angle=270,width=12.2cm]{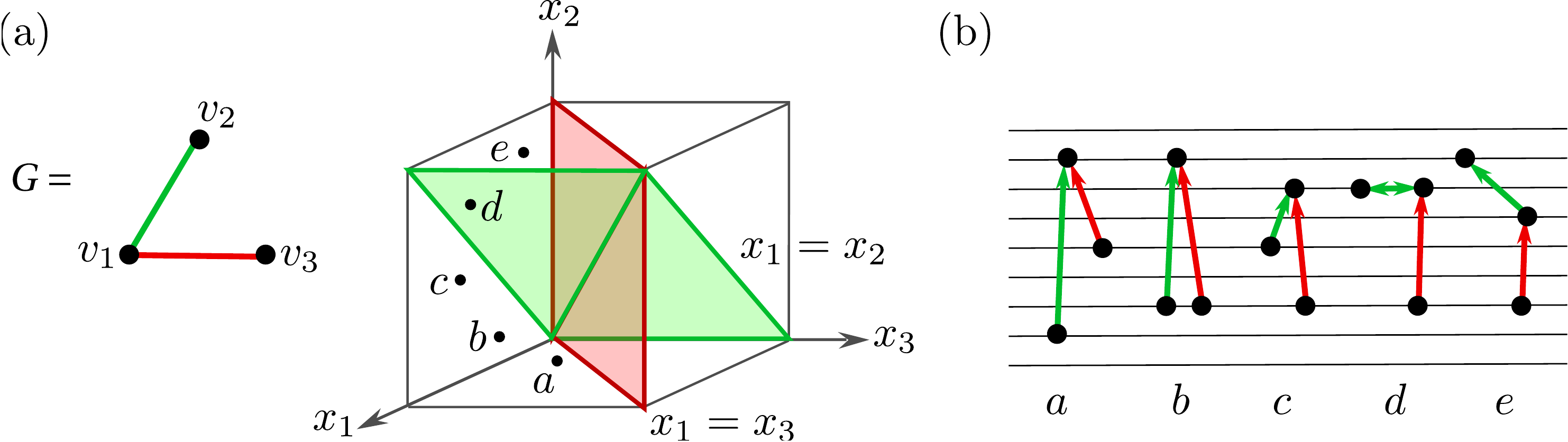}
\caption{(a) The graphic arrangement of the graph $G$. An edge between vertices $v_i$ and $v_j$ corresponds to a hyperplane $x_i=x_j$. (b) Points in the cube correspond to colorings. They can be visualized by drawing the graph $G$ in a coordinate system such that the height of a vertex $v_i$ is given by its color $x_i\in\{1,\ldots,k\}$. This induces an orientation of the edges from the vertex with smaller to the vertex with the larger color. Moving from coloring $a$ through coloring $b$ to coloring $c$ we pass the hyperplane $x_2=x_3$ which is not part of the graphic arrangement; in $b$ vertices $v_2$ and $v_3$ are at the same height. Moving on through $d$ to $e$ we pass the hyperplane $x_1=x_2$ which is in the graphic arrangement; in $d$ two adjacent vertices are at the same height, so $d$ is not proper. Moving from $c$ to $e$ thus reverses the orientation of the edge between vertices $v_1$ and $v_2$.}
\label{fig:chromatic}
\end{figure}

The chromatic polynomial $\chi_G(k)$ of a given graph $G$ counts the number of proper $k$-colorings of $G$. Let $V$ be the vertex set of $G$ and $\sim$ its adjacency relation. A \emph{$k$-coloring} is a vector $x\in [k]^V$ that assigns to each vertex $v\in V$ a color $x_v\in [k] := \{1,\ldots,k\}$. Such a $k$-coloring $x$ is \emph{proper} if for any two adjacent vertices $v\sim w$ the assigned colors are different, i.e., $x_v \not= x_w$. This way of describing a coloring as a vector rather than a function already suggests a geometric point of view (Figure~\ref{fig:chromatic}). Define the \emph{graphic arrangement} of $G$ as the set of all hyperplanes $x_v = x_w$ for adjacent vertices $v \sim w$. Then the chromatic polynomial counts integer points $x\in\ZZ^V$ that are contained in the half-open cube $(0,k]^V$ but do not lie on any of the hyperplanes in the graphic arrangement of $G$, i.e.,
\begin{eqnarray}
\label{eqn:chromatic}
\chi_G(k) = \#\ZZ^V\cap\mset{x\in\RR^V}{0< x_v \leq k \text{ and }  x_v \not= x_w \text{ if } v\sim w}.
\end{eqnarray}

\begin{figure}[t]
\includegraphics[angle=270,width=12.2cm]{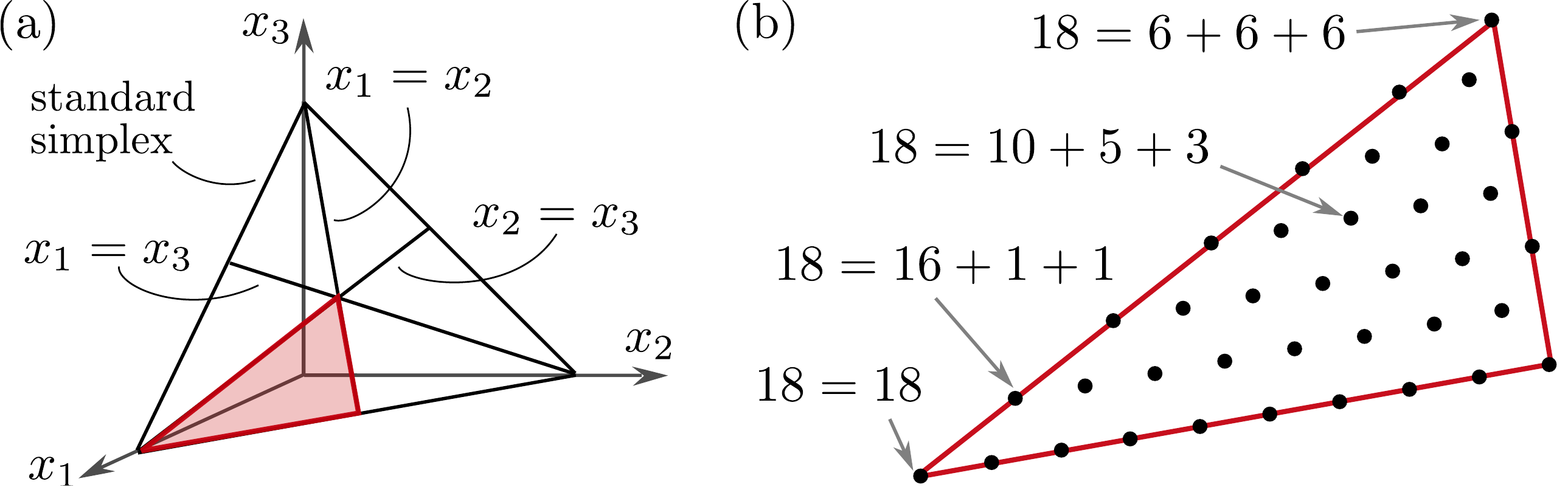}
\caption{(a) The partition polytope is cut out from the standard simplex $\mset{x}{x_i\geq 0, \sum x_i = 1}$ by the braid arrangement of all hyperplanes $x_i=x_j$. (b) The integer points in the partition polytope for $k=18$ and $m=3$ correspond to the partitions of 18 in at most 3 parts.}
\label{fig:partition}
\end{figure}

The restricted partition function $p(k,m)$ counts the number of partitions of $k$ into at most $m$ parts.\footnote{It is easy to adapt the following construction to the case of counting partitions with exactly $m$ parts by making one inequality strict.} This can be modeled simply by defining a partition of $k$ into at most $m$ parts as a non-negative vector $x\in\ZZ^m$ whose entries sum to $k$ and are weakly decreasing. For example, in the case $m=5$ and $k=14$ the partition $14=7+5+2$ would correspond to the vector $(7,5,2,0,0)$. In short,
\begin{eqnarray}
\label{eqn:partition}
p(k,m) = \#\ZZ^m\cap\mmset{x\in\RR^m}{x_1 \geq x_2 \geq \ldots \geq x_m \geq 0 \text{ and }  \sum_{i=1}^m x_i = k}.
\end{eqnarray}
Geometrically speaking, the restricted partition function thus counts integer points in an $(m-1)$-dimensional simplex in $m$-dimensional space. This is visualized in Figure~\ref{fig:partition}. Note that the constraints that all variables are non-negative and that their sum is equal to $k$ already defines an $(m-1)$-dimensional simplex, bounded by the coordinate hyperplanes. The \emph{braid arrangement}, i.e., the set of all hyperplanes $x_i=x_j$, subdivides this simplex into $m!$ equivalent pieces; the definition of the restricted partition function then selects the one piece in which the coordinates are in weakly decreasing order.

It is interesting to observe that both constructions work with the braid arrangement. Indeed, there are a host of combinatorial models that fit into this setting. A great example are \emph{scheduling problems} \cite{Breuer2014}: Given a number $k$ of time-slots, how many ways are there to schedule $d$ jobs such that they satisfy a boolean formula $\psi$ over the atomic expressions ``job $i$ runs before job $j$'', i.e., $x_i < x_j$?  E.g., if $\psi=(x_1< x_2)\rar(x_3< x_2)$ then we would count all ways to place $3$ jobs in $k$ time-slots such that if job $1$ runs before job $2$, then job $3$ also has to run before $2$. We will return to scheduling problems in Section~\ref{sec:quasisymmetric}. However, the methods presented in this article are not restricted to this setup as we will see.

\section{Ehrhart Theory}
\label{sec:ehrhart}

For any set $X\subset\RR^d$ the \emph{Ehrhart function} $\ehr_X(k)$ of $X$ counts the number of integer points in the $k$-th dilate of $X$ for each $1\leq k\in \ZZ$, i.e.,
\[
	\ehr_X(k) = \#\ZZ^d \cap (k\cdot X).
\]
Both our constructions from the previous section are of this form, since (\ref{eqn:chromatic}) and (\ref{eqn:partition}) are, respectively, equivalent to
\begin{eqnarray*}
\chi_G(k) & =& \#\ZZ^V\cap k\cdot\mset{x\in\RR^n}{0< x_v \leq 1 \text{ and }  x_v \not= x_w \text{ if } v\sim w}, \\
p(k,m) &=& \#\ZZ^m\cap k\cdot\mmset{x\in\RR^m}{x_1 \geq x_2 \geq \ldots \geq x_m \geq 0 \text{ and }  \sum_{i=1}^m x_i = 1}.
\end{eqnarray*}

We will call the set $X$ the \emph{geometric model} of the counting function $\ehr_X$. The central theme of this exposition is that geometric properties of $X$ often translate into algebraic properties of $\ehr_X$. Ehrhart's theorem is the prime example of this phenomenon. To set the stage, we introduce some terminology and refer to \cite{Schrijver1986,Ziegler1995} for concepts from polyhedral geometry not defined here.

A \emph{polyhedron} is any set of the form $P=\mset{x \in\RR^d}{Ax \geq b}$ for a fixed matrix $A$ and vector $b$. All polyhedra in this article will be \emph{rational}, i.e., we can assume that $A$ and $b$ have only integral entries. A \emph{polytope} is a bounded polyhedron. Any dilate of a polytope contains only a finite number of integer points, whence the Ehrhart function of a polytope is well-defined. A polyhedron is \emph{half-open} if some of its defining inequalities are strict. A \emph{partial polyhedral complex}\footnote{Classically, a \emph{polyhedral complex} is a collection $X$ of polyhedra that is closed under passing to faces, such that the intersection of any two polyhedra in $X$ is also in $X$ and is a face of both. In contrast, in a partial polyhedral complex some faces are allowed to be open. This means that it is possible to remove an edge from a triangle -- including or excluding the incident vertices. It is sometimes useful to regard a partial polyhedral complex as subset of a fixed underlying polyhedral complex, so as to be able to refer to the vertices of the underlying complex, for example. We will disregard these technical issues in this expository paper, however.} $X$ is any set that can be written as a disjoint union of half-open polytopes.

Our model of $p(k,m)$ is a polytope. Our model of $\chi_G(k)$ is not, though, as it is non-convex, disconnected and neither closed nor open. It is easily seen to be a partial polytopal complex, though, e.g., by rewriting $x_v \not= x_w$ to $(x_v < x_w) \vee (x_v > x_w)$ and bringing the resulting formula in disjunctive normal form. This makes partial polytopal complexes an extremely flexible modeling framework, as summarized in the following lemma.

\begin{lemma}
\label{lem:partial-polytopal-complex}
Let $\psi$ be any boolean formula over homogeneous linear equations and inequalities with rational coefficients in the variables $x_1,\ldots,x_d,k$, such that for every $k$ the set of all $x$ such that $\psi(x)$ is bounded. Then there exists a partial polytopal complex $X$ such that for all $1 \leq k\in\ZZ$,
\[
	\#\mset{x\in\ZZ^d}{\psi(x,k)} = \ehr_X(k).
\]
\end{lemma}

The generality of Ehrhart functions of partial polytopal complexes underlines the strength of the following famous theorem by Eug\`ene Ehrhart.

\begin{theorem}[Ehrhart~\cite{Ehrhart1962}]
\label{thm:ehrhart}
If $X$ is partial polytopal complex\footnote{Ehrhart formulated his theorem for polytopes, not for partial polytopal complexes. The generalization follows immediately, however, since for any partial polytopal complex $X$ the Ehrhart function $\ehr_X$ is a linear combination of Ehrhart functions of polytopes.}, then $\ehr_X(k)$ is a quasipolynomial.
\end{theorem} 

\begin{figure}[t]
\includegraphics[angle=270,width=12.2cm]{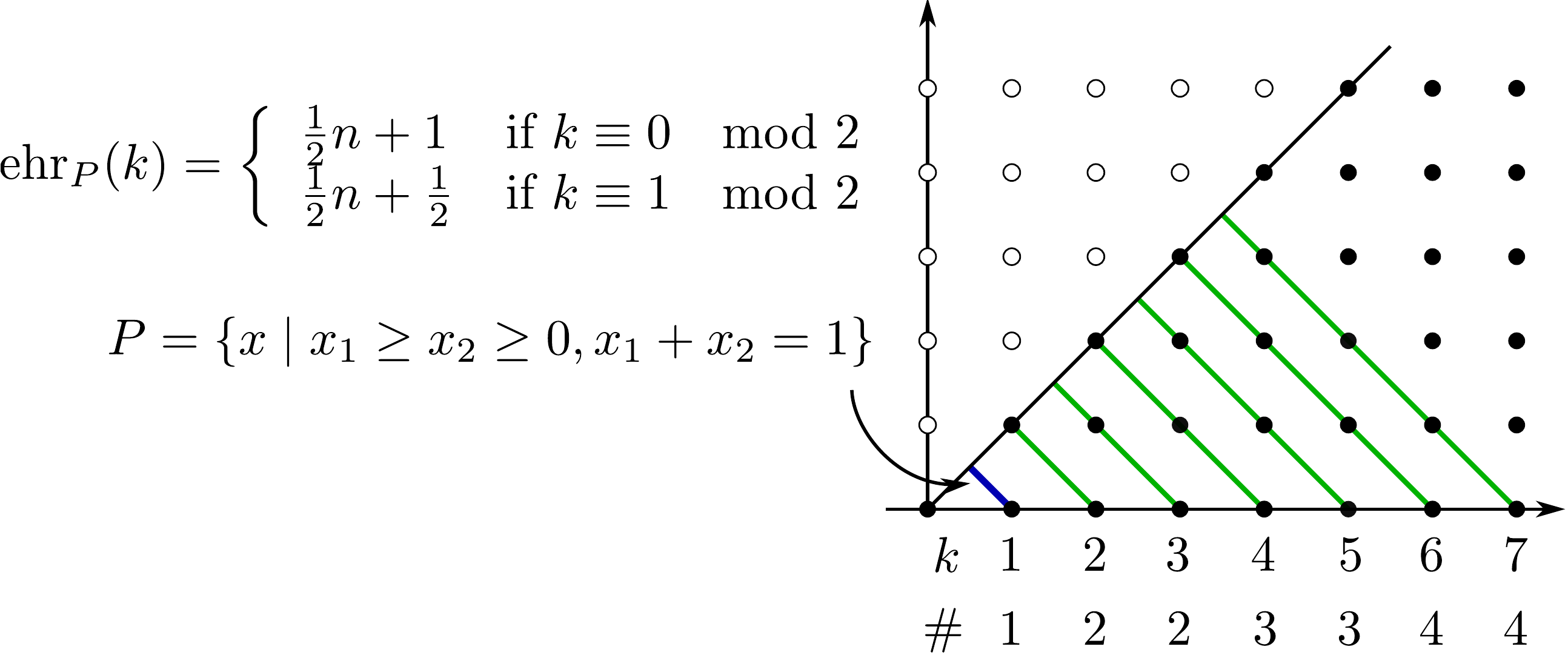}
\caption{By counting the integer points in the dilates of $P$ and interpolating, we can compute the Ehrhart quasipolynomial of $P$. In this case $P$ is the restricted partition polytope for partitions into at most 2 parts.}
\label{fig:quasipolynomial}
\end{figure}

Quasipolynomials are an important class of counting functions which capture both polynomial growth and periodic behavior. A function $p(k)$ is a \emph{quasipolynomial} if there exist polynomials $p_0(k),\ldots,p_{\ell-1}(k)$ such that
\[
	p(k) = 
		\choice{
			p_0(k) & \text{if } k \equiv 0 \mod \ell \\
			p_1(k) & \text{if } k \equiv 1 \mod \ell \\
			\vdots \\
			p_{\ell-1}(k) & \text{if } k \equiv \ell-1 \mod \ell
		}
\]
for all $k\in\ZZ$. The polynomials $p_i$ are called the \emph{constituents} of $p$ and their number is a \emph{period} of $p$. The period is not uniquely determined, but of course the \emph{minimal period} is; every period is a multiple of the minimal period. The \emph{degree} of $p$ is the maximal degree of the $p_i$. If $d$ is the degree of $p$ and $\ell$ a period of $p$, then $p$ is uniquely determined by $\ell\cdot (d+1)$ values of $p$, or more precisely, by $d+1$ values $p(k'\cdot\ell + i)$ for each $i=0,\ldots,\ell-1$. This is why it makes sense to say that $\ehr_X$ ``is'' a quasipolynomial, even though we have defined Ehrhart functions only at positive integers. As an example, the quasipolynomial given by the restricted partition function $p(k,2)$ is computed by interpolation in Figure~\ref{fig:quasipolynomial}.

Given this terminology, we can make our above statement of Ehrhart's theorem more precise. Restricting our attention to polytopes $P$ for the moment, the following hold for $\ehr_P$. First, all constituents of $\ehr_P$ have the same degree. That degree is the dimension of $P$. Second, the leading coefficient of all constituents of $P$ in the monomial basis is the volume of $P$. Third, the least common multiple of the denominators of all vertices of $P$ is a period of $\ehr_P$. More precisely, if $v_1,\ldots,v_N\in\QQ^d$ are the vertices of $P$ and $v_{i,j}=\frac{a_{i,j}}{b_{i,j}}\in\QQ$, then
\[
\ell=\lcm(\mset{b_{i,j}}{i=1,\ldots,N, \; j=1,\ldots,d})
\]
is a period of $\ehr_P$. In particular, if the vertices of $P$ are all integral the Ehrhart function is a polynomial.

Ehrhart's theorem thus provides a very general method for proving that counting functions are (quasi-)polynomials. Simply by virtue of the geometric models from Section~\ref{sec:modeling}, we immediately obtain that the chromatic function is a polynomial because the vertices of our geometric model are integral. This proof of polynomiality is very different from the standard deletion-contraction method and generalizes to counting functions that do not satisfy such a recurrence, including all scheduling problems. Also we find that the restricted partition function into $m$ parts is a quasipolynomial with period $\lcm(1,2,\ldots,m)$: The numbers $1,2,\ldots,m$ appear in the denominators of the vertices, since we intersect the braid arrangement with the simplex $\mset{x}{x_i\geq 0,\sum x_i =1}$ instead of the cube. For more on the restricted partition function from an Ehrhart perspective, see \cite{BEKZ14}.

Lemma~\ref{lem:partial-polytopal-complex} can be generalized even further, for example by allowing quantifiers via Presburger arithmetic \cite{Woods2013} or by considering the multivariate case \cite{Verdoolaege2007}. As a great introductory textbook on Ehrhart theory we recommend \cite{Beck2007}.

\section{Combinatorial Reciprocity Theorems}
\label{sec:reciprocity}

Now that we know that the Ehrhart function $\ehr_P$ of a polytope $P$ is in fact a quasipolynomial we can evaluate it at negative integers. Even though the Ehrhart function itself is defined only at positive integers, it turns out that the values of $\ehr_P$ at negative integers have a very elegant geometric interpretation: $\ehr_P(-k)$ counts the number of integer points in the interior of $k\cdot P$.

\begin{theorem}[Ehrhart-Macdonald Reciprocity \cite{Macdonald1971}]
\label{thm:reciprocity}
If $P\subset\RR^n$ is a polytope of dimension $d$ and $0<k \in\ZZ$ then
\[
	\ehr_P(-k) = (-1)^d \cdot (\#\ZZ^n \cap k\cdot P^\circ).
\]
\end{theorem}

Here $P^\circ$ denotes the \emph{relative interior} of $P$, which means the interior of $P$ taken with respect to the affine hull\footnote{The affine hull of $P$ is the smallest affine space containing $P$. Affine spaces are the translates of linear spaces.} of $P$. If $P$ is given in terms of a system of linear equations and inequalities, the relative interior is often easy to determine. For example, if $P$ is defined by $Ax \geq b$ and $A'x = b'$, and $A'$ contains all equalities of the system\footnote{More precisely, we require that the affine hull of $P$ is $\mset{x}{A'x=b'}$ and that for every row $a$ of $A$ the linear functional $\sprod{a}{x}$ is not constant over $x\in P$.}, then $P^\circ$ is given by $Ax > b$ and $A'x = b'$. In short, all we need to do is make weak inequalities strict.

Ehrhart-Macdonald reciprocity provides us with a powerful framework for finding combinatorial reciprocity theorems, i.e., combinatorial interpretations of the values of counting functions at negative integers. We start with a counting function $f$ defined in the language of combinatorics and translate this counting function into the language of geometry by constructing a linear model. In the world of geometry, we apply Ehrhart-Macdonald reciprocity to find a geometric interpretation of the values of $f$ at negative integers. Translating this geometric interpretation back into the language of combinatorics, a process which can be quite subtle, we then arrive at a combinatorial reciprocity theorem.

Let us start with the example of the restricted partition function $p(k,m)$. Applying Theorem~\ref{thm:reciprocity} it follows that, up to sign, $p(-k,m)$ counts vectors $x$ such that $x_1 > x_2 > \ldots > x_m > 0$ and $\sum x_i = k$ for any positive integer $k$. Interpreting this geometric statement combinatorially, we find:

\begin{theorem}
Up to sign, $p(-k,m)$ counts partitions of $k$ into exactly $m$ distinct parts.
\end{theorem}

This result seems to be less well-known in partition theory than one would expect, even though it is an immediate consequence of Ehrhart-Macdonald reciprocity; see also \cite{BEKZ14}. A very similar geometric construction, however, is the basis of Stanley's work on $P$-partitions and the order polynomial \cite{Stanley1986} which has many nice extensions, e.g., \cite{Jochemko2013}.

Next, we consider the chromatic polynomial $\chi_G(k)$. The model $X_G$ we use here is slightly different from (\ref{eqn:chromatic}) in that we work with the open cube $(0,k+1)^V$. This introduces a shift $\ehr_{X_G}(k) = \chi_G(k-1)$. The advantage is that $X_G$ is now a disjoint union of open polytopes $P_1,\ldots,P_N$. As already motivated by Figure~\ref{fig:chromatic}, it turns out that the $P_i$ are in one-to-one correspondence with the acyclic orientations\footnote{An orientation of a graph $G$ is \emph{acyclic}, if it contains no directed cycles.} of the graph $G$ \cite{Greene1977}. Applying Theorem~\ref{thm:reciprocity} to each component individually, we find that $\chi_G(-k)$ counts all integer vectors $x$ in the closed cube such that points on the hyperplanes $x_v=x_w$ have a multiplicity equal to the number of closed components $\bar{P}_i$ they are contained in. To interpret this combinatorially, we define an orientation $o$ and a coloring $x$ of $G$ to be \emph{compatible} if, when moving along directed edges, the colors of the vertices always increase or stay the same. Putting everything together and taking the shift into account we obtain Stanley's reciprocity theorem for the chromatic polynomial below. The geometric proof we described is due to Beck and Zaslavsky \cite{Beck2006-iop} and can be generalized to cell-complexes \cite{BBGM2014}.

% ehr(k) = chi(k-1) => ehr(-k) = chi(-k-1) => chi(-k) = ehr(-k+1) <-> [0,k-1] closed <-> k colors

\begin{theorem}[Stanley \cite{Stanley1973}]
Up to sign, $\chi_G(-k)$ counts pairs $(x,o)$ of (not necessarily proper) $k$-colorings $x$ and compatible acyclic orientations $o$ of $G$.
\end{theorem}

\begin{figure}[t]
\includegraphics[angle=270,width=12.2cm]{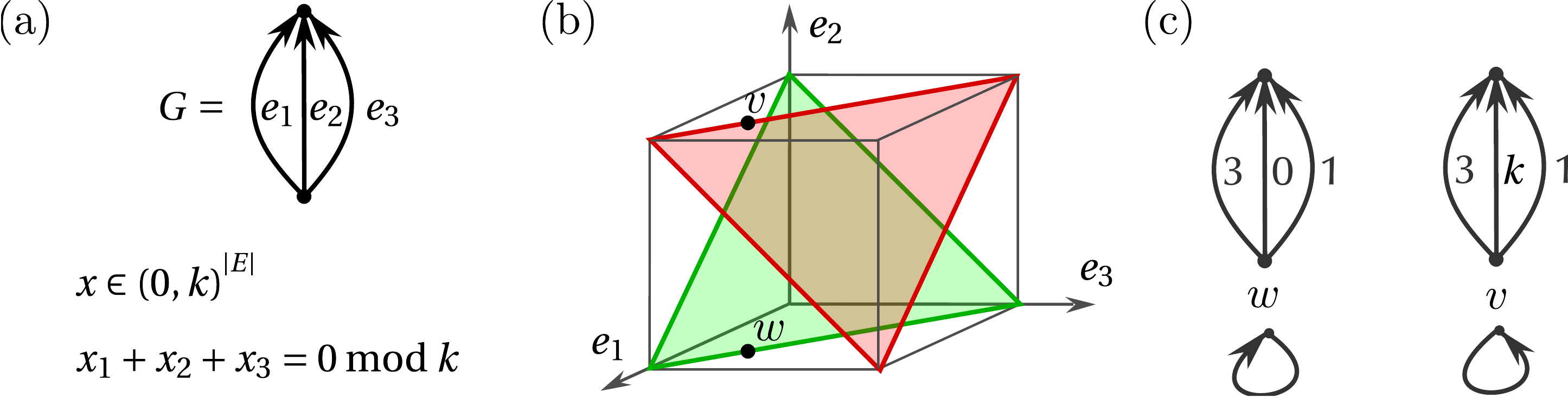}
\caption{(a) A directed graph $G$ and the flow problem $G$ defines. (b) The corresponding partial polyhedral complex $X$. (c) The labelings of $G$ given by the points $v,w$ along with the different totally cyclic orientations of $G/\supp(v) = G/\supp(w)$.}
\label{fig:flow}
\end{figure}

Next, the modular flow polynomial of a graph provides us with an example of a combinatorial reciprocity theorem that was first discovered via the geometric approach and that makes use of a different construction, unrelated to the braid arrangement. This example is illustrated in Figure~\ref{fig:flow}. A $\ZZ_k$-flow on a directed graph $G$ with edge set $E$ is a vector $y\in\ZZ_k^E$ that assigns to each edge of $G$ a number such that at each vertex $v$ of $G$ the sum of all flows into $v$ equals the sum of all flows out of $v$, modulo $k$. The modular flow polynomial $\varphi_G(k)$ of $G$ counts $\ZZ_k$-flows on $G$ that are nowhere zero. To model this in Euclidean space, we identify the elements of $\ZZ_k$ with the integers $0,\ldots,k-1$. Nowhere zero vectors $y\in\ZZ_k^E$ thus correspond to integer points $y\in(0,k)^k$ in the $k$-th dilate of the open unit cube. If $A\in\ZZ^{V\times E}$ is the incidence matrix of $G$, the constraint that flow has to be conserved at each vertex can be expressed simply by requiring $Ay\equiv 0 \mod k$, or, equivalently, by $\exists b\in\ZZ^V: Ay = kb$. Note that for only finitely many $b\in\ZZ^V$ the hyperplane $Ay=b$ intersects the unit cube $(0,1)^E$. Let $P_1,\ldots,P_N$ denote these sections and let $X$ be their union. Then $\varphi_G(k) = \ehr_X(k)$.

Applying Ehrhart-Macdonald reciprocity, we obtain that, up to sign, $\varphi_G(-k)$ counts integer points in the $k$-th dilate of the union of the closures $\bar{P}_i$. In particular, we now count vectors that may have both entries $0$ and $k$, which are both congruent zero mod $k$, but which we have to count as different as Figure~\ref{fig:flow} shows. This observation suggests that to find a combinatorial interpretation, we may want to consider assigning two different kinds of labels to the edges with zero flow. Pursuing this line of thought eventually leads to the following combinatorial reciprocity theorem, which again can be generalized to cell complexes \cite{BBGM2014}, see also \cite{Beck2006-flow,Breuer2009}.

\begin{theorem}[Breuer-Sanyal \cite{BreuerSanyal2012}]
Up to sign, $\varphi_G(-k)$ counts pairs $(y,o)$ of a $\ZZ_k$-flow $y$ on $G$ and a totally cyclic reorientation of $G/\supp(y)$.
\end{theorem}

Here a reorientation is a labeling of the edges of a directed graph with $+$ or $-$, indicating whether the direction of the edge should be reversed or not. Such a reorientation is \emph{totally cyclic}, if every edge of the resulting directed graph lies on a directed cycle. $\supp(y)$ denotes the set of edges where $y$ is non-zero and $G/\supp(y)$ denotes the graph where $\supp(y)$ has been contracted. 

For more on combinatorial reciprocity theorems we recommend the forthcoming book~\cite{Beck2014}.

\section{Cones and Fundamental Parallelepipeds}
\label{sec:cones}

A polyhedral cone or \emph{cone}, for short, is the set of all linear combinations with non-negative real coefficients of a finite set of \emph{generators} $v_1,\ldots,v_d\in\QQ^n$. If the generators are linearly independent, the cone is \emph{simplicial}. The cone is \emph{pointed} or \emph{line-free} if it does not contain a line $\{u+\lambda v\;|\; \lambda \in\RR\}$.

\begin{figure}[t]
\includegraphics[angle=270,width=12.2cm]{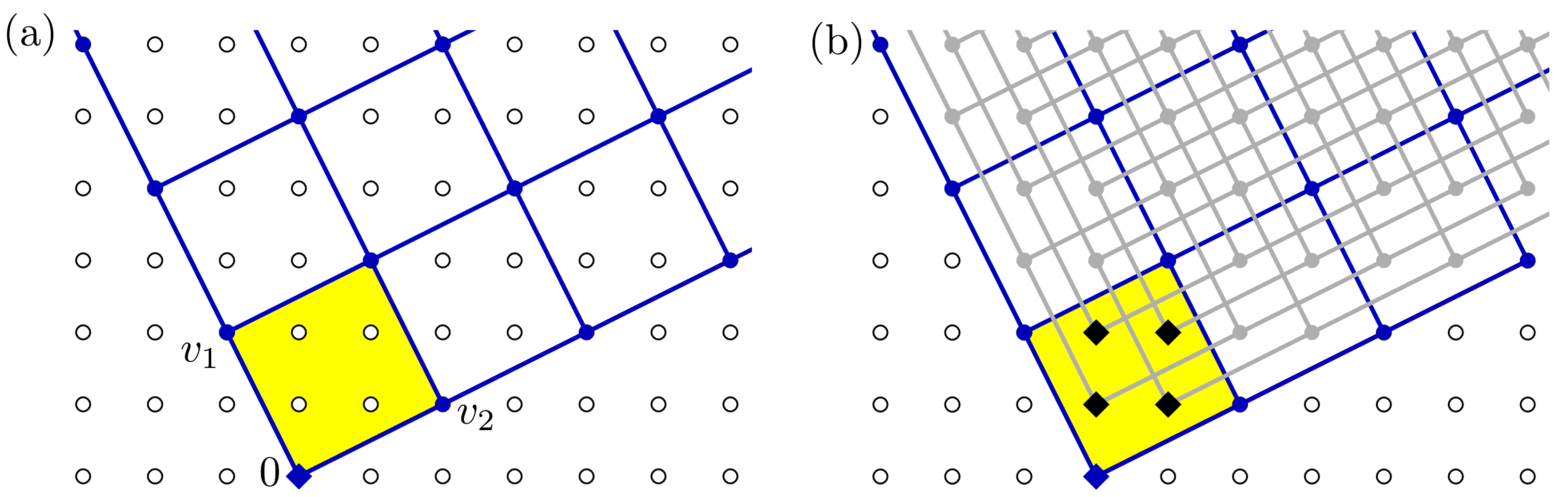}
\caption{(a) The discrete cone $\cone_\ZZ(v_1,v_2)$ generated by $v_1$ and $v_2$. Its fundamental parallelepiped is shaded. (b) To generate all points in $\ZZ^2\cap\cone_\RR(v_1,v_2)$ the discrete cone has to be translated by every integer point in the fundamental parallelepiped (shown as diamonds).}
\label{fig:fundamental-lemma}
\end{figure}

Cones are the basic building blocks of Ehrhart theory, because the sets of integer points in simplicial cones have a very elegant description, which is illustrated in Figure~\ref{fig:fundamental-lemma}. Let $v_1,\ldots,v_d\in\ZZ^n$ be linearly independent, and consider the simplicial cone $\cone_\RR(v_1,\ldots,v_d)$ generated by them. The \emph{discrete cone} or \emph{semigroup} $\cone_\ZZ(v_1,\ldots,v_d)$ of all non-negative \emph{integral} combinations of the $v_i$ reaches only those integer points in $\cone_\RR$ that lie on the lattice $\ZZ v_1+\ldots + \ZZ v_d$ generated by the $v_i$. However, by shifting the discrete cone to all integer points in the \emph{fundamental parallelepiped} $\Pi(v_1,\ldots,v_d)$
we can not only capture all integer points in $C$, but we moreover partition them into $\#\ZZ^{n+1}\cap\Pi(v_1,\ldots,v_d) = |\det(v_1,\ldots,v_d)|$ disjoint classes. This number of integer points in the fundamental parallelepiped is called the \emph{index} of $C$. Define
\begin{eqnarray*}
\cone_\RR(v_1,\ldots,v_d) & = & \mset{\sum_{i=1}^d \lambda_i v_i}{0\leq \lambda_i \in \RR}, \\
\cone_\ZZ(v_1,\ldots,v_d) & = & \mset{\sum_{i=1}^d \lambda_i v_i}{0\leq \lambda_i \in \ZZ}, \\
\Pi(v_1,\ldots,v_d) & = & \mset{\sum_{i=1}^d \lambda_i v_i}{0\leq \lambda_i < 1}.
\end{eqnarray*}

\begin{lemma}
\label{lem:fundamental-lemma}
Let $v_1,\ldots,v_d\in\ZZ^n$ be linearly independent. Then 
\[
  \ZZ^n\cap\cone_\RR(v_1,\ldots,v_d) = \left(\ZZ^n\cap\Pi(v_1,\ldots,v_d)\right) + \cone_\ZZ(v_1,\ldots,v_d).
\]
\end{lemma}

The main benefit of this decomposition is that it splits the problem of describing the integer points in a cone to into two parts: The finite problem of enumerating the integer points in the fundamental parallelepiped, and the problem of describing the discrete cone, which is easy as we shall see below. As an application of this result, we will now prove Ehrhart's theorem for polytopes.

\begin{figure}[t]
\begin{center}
\includegraphics[angle=270,width=10cm]{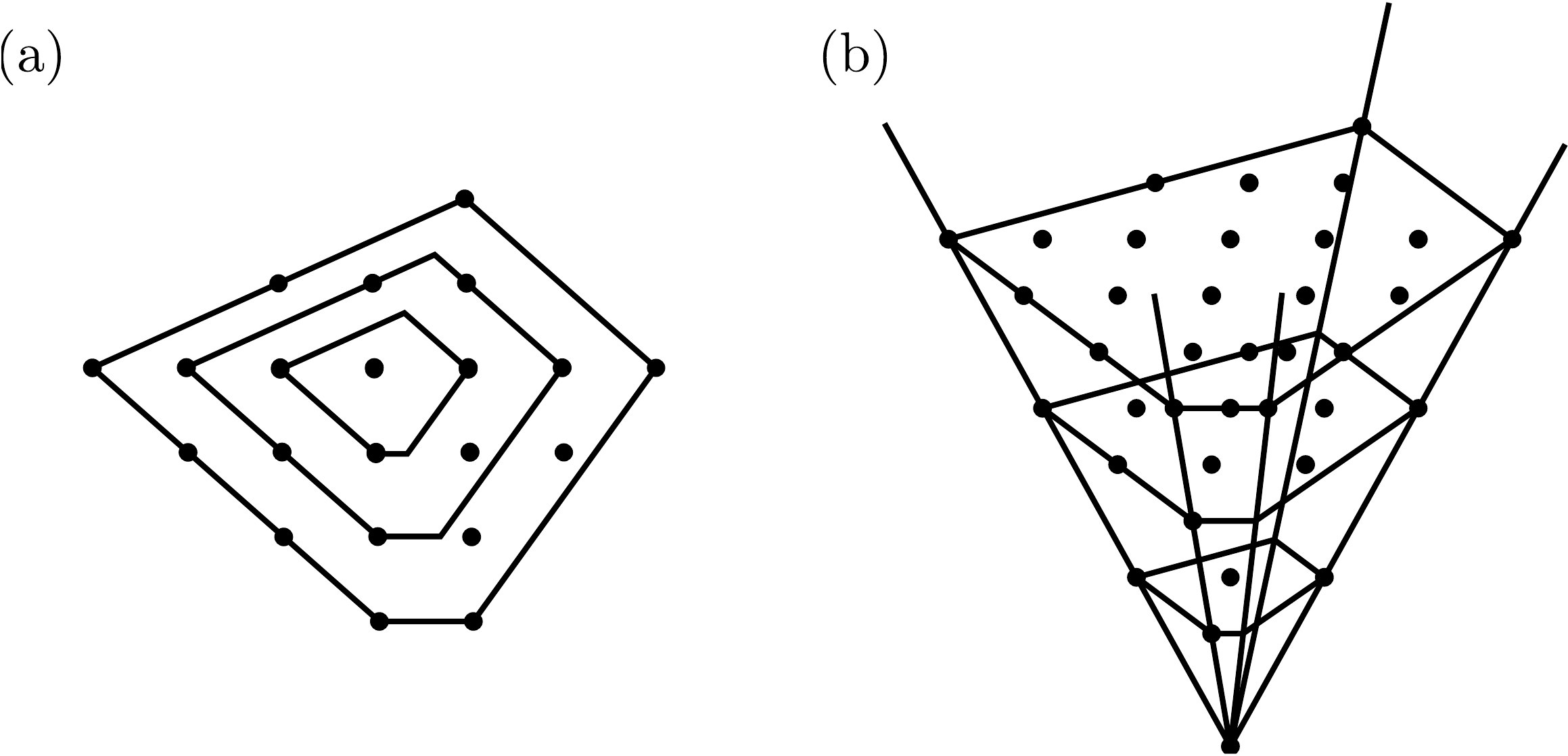}
\end{center}
\caption{(a) The three first dilates of a polytope $P'$ in the plane. (b) The polytope $P=P'\times\{1\}$ embedded in 3-space and the cone $C=\cone(P)$ over $P$. Sections $H_k\cap C$ are lattice equivalent to the dilates of $P'$.}
\label{fig:ehrhart-cone}
\end{figure}

Suppose we want to compute the Ehrhart function of a polytope $P'\subset\RR^n$. We embed $P'$ at height $1$ in $\RR^{n+1}$, i.e., we pass to $P = P'\times\{1\}\subset\RR^{n+1}$. Then, we consider the set $\cone(P)$ of all finite linear combinations of elements in $P$ with non-negative real coefficients as shown in Figure~\ref{fig:ehrhart-cone}. The intersections of $\cone(P)$ with the hyperplanes $H_k:=\mset{x}{x_{n+1} = k}$ are lattice equivalent\footnote{Two sets $X,Y\subset\ZZ^n$ are \emph{lattice equivalent} if there exists an affine isomorphism $x \mapsto Ax+b$ that maps $X$ to $Y$ and which induces a bijection on $\ZZ^n$.} to the dilates $k\cdot P$ we are interested in. If we can describe the number of integer points in such sections of polyhedral cones, we will have a handle on computing Ehrhart functions. 

Before we continue, we observe that we can make two more simplifications. First, we can restrict our attention to simplicial cones. While in general $\cone(P)$ will of course not be simplicial, we can always reduce the problem to simplicial cones by triangulating $P$. Second, we note that while $\cone(P)$ is indeed finitely generated by the vertices $w_1,\ldots,w_N$ of $P$, these $w_i$ may be rational vectors. Instead, we would like to work with generators $v_i$ that are all integer and all at the same height wrt.~the last coordinate. This can be achieved by letting $\ell$ denote the smallest integer such that $\ell\cdot w_i\in\ZZ^{n+1}$ for all $i$ and setting $v_i := \ell\cdot w_i$.

\begin{figure}[t]
\begin{center}
\includegraphics[angle=270,width=10cm]{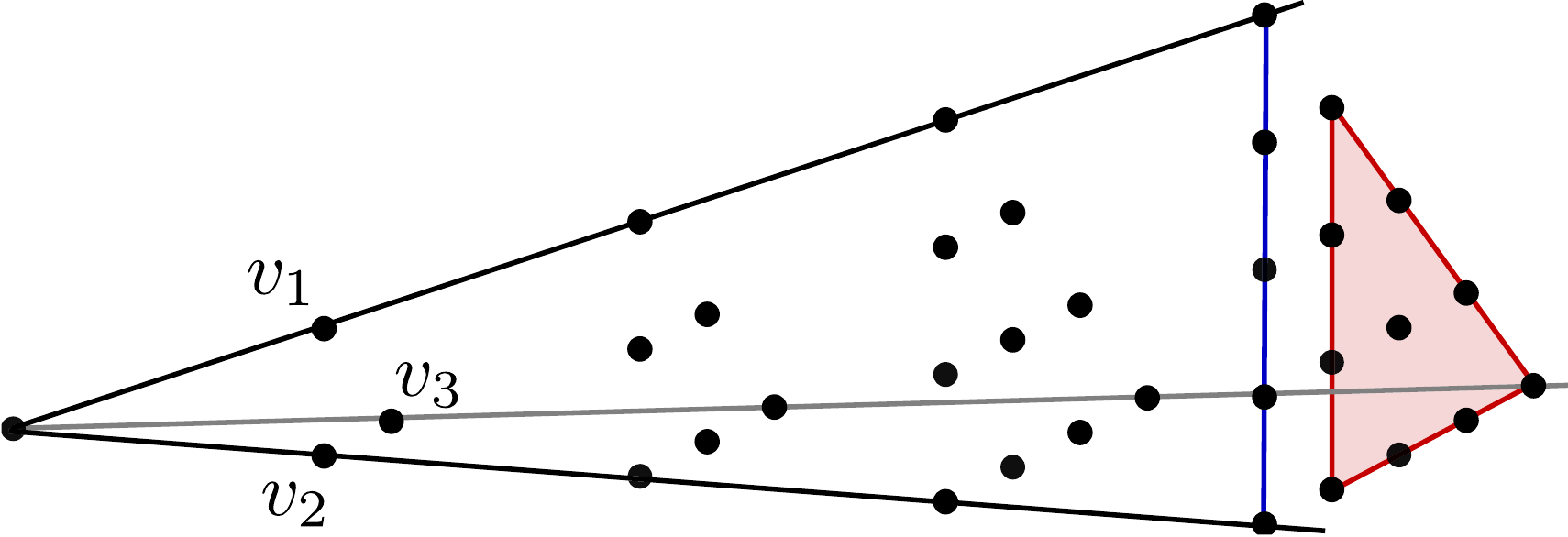}
\end{center}
\caption{The $4\ell$-th level of $C=\cone_\ZZ(v_1,v_2,v_3)$ decomposes naturally into the $4\ell$-th level of $\cone_\ZZ(v_1,v_2)$ and a shift of the $3\ell$-th level of $C$.}
\label{fig:pascal}
\end{figure}

We have thus reduced the problem of computing the Ehrhart function of $P'$ to computing $\ZZ^{n+1}\cap H_k\cap C_\RR$, the number of integer points at height $k$ in a simplicial cone $C_\RR:=\cone_\RR(v_1,\ldots,v_d)$ given by integral generators with last coordinate equal to a constant $\ell$. Following Lemma~\ref{lem:fundamental-lemma} we concentrate on $C_\ZZ:=\cone_\ZZ(v_1,\ldots,v_d)$ first. Since the last coordinate of all $v_i$ is $\ell$, $C_\ZZ\cap H_k$ is empty if $k \not\equiv 0 \mod \ell$. On the other hand, if $k \equiv 0 \mod \ell$ and $k>0$, then
\[
	H_{k} \cap \cone_\ZZ(v_1,\ldots,v_d) = v_d + (H_{k-\ell} \cap \cone_\ZZ(v_1,\ldots,v_d)) \cup (H_{k} \cap \cone_\ZZ(v_1,\ldots,v_{d-1})).
\]
This is an instance of Pascal's recurrence for the binomial coefficients, illustrated in Figure~\ref{fig:pascal}, which yields for all integers $k\geq 0$,
\[
  \#H_{\ell \cdot k} \cap \cone_\ZZ(v_1,\ldots,v_d) = \binom{k+d-1}{d-1}
  %\choice{ 
  %  \binom{k+d-1}{d-1} & \text{if } k\equiv 0 \mod \ell, \\ 
  %  0 & \text{if } k\not\equiv 0 \mod \ell. 
  %}
\]
and $\#H_{k} \cap \cone_\ZZ(v_1,\ldots,v_d) = 0$ if $k \not\equiv 0 \mmod \ell$.

Applying Lemma~\ref{lem:fundamental-lemma} we see that to get the counting function for $C_\RR$, we need to shift the discrete cone $C_\ZZ$ by all the integer points in the fundamental parallelepiped, which allows us to reach lattice points at heights $k$ which are not a multiple of $\ell$. Organizing these shifts according to the last coordinate, we obtain for any $k\geq 0$ and $0\leq r < \ell$
\begin{eqnarray}
	&& \#\left(\ZZ^{n+1}\cap H_{\ell\cdot k+r} \cap C_\RR\right) \nonumber \\
  &=& h^*_r \binom{k+d-1}{d-1} + h^*_{\ell+r} \binom{k+d-2}{d-1} + \ldots + h^*_{(d-1)\cdot\ell+r} \binom{k}{d-1}
	\label{eqn:counting-slices}
\end{eqnarray}
where $h^*_i$ denotes the number of integer points at height $i$ in $\Pi(v_1,\ldots,v_d)$. 

Note that if we are interested in $\#\left(\ZZ^{n+1}\cap H_{m} \cap C_\RR\right)$ for an arbitrary non-negative $m$ then we can always write $m = \ell k + r$ such that $k\geq 0$ and $0\leq r < \ell$ simply by doing division with remainder. Also note that (\ref{eqn:counting-slices}) is a polynomial of degree $d-1$ in $k$ for each fixed $r$. Since $r$ changes periodically with $m$, the counting function $m \mapsto \#\left(\ZZ^{n+1}\cap H_{m} \cap C_\RR\right)$ is a quasipolynomial of period $\ell$. By construction, $\ehr_P(k)$ is a sum of such expressions and therefore itself a quasipolynomial, which completes the proof of Theorem~\ref{thm:ehrhart}.

\section{Connection to Rational Functions}
\label{sec:brion}

The results and constructions of the previous section translate immediately into the language of generating functions, formal power series and rational functions. When we represent an integer point $v\in\ZZ^n$ by a multivariate monomial $z^v:=z_1^{v_1}\cdot\ldots\cdot z_n^{v_n}$, the set of integer vectors in any given set $S\subset\ZZ^n$ can be written as a multivariate generating function
\[
	\phi_S(z) = \sum_{v\in \ZZ^n \cap S} z^v.
\]
Using the familiar geometric series expansion $\frac{1}{1-z^v} = \sum_{i=0}^\infty z^{iv}$ we see that generating functions of ``discrete rays'' of integer vectors can be represented as rational functions. Indeed, both discrete cones and Lemma~\ref{lem:fundamental-lemma} can be expressed succinctly in terms of rational functions.
\begin{eqnarray}
	\phi_{\cone_\ZZ(v_1,\ldots,v_d)}(z)
	& = &
	\frac
		{1}
		{(1-z^{v_1})\cdot\ldots\cdot(1-z^{v_d})}.
	\label{eqn:rat-fun-disc-cone}
\\
	\phi_{\cone_\RR(v_1,\ldots,v_d)}(z)
	& = &
	\frac
		{\sum_{v\in \ZZ^n\cap\Pi(v_1,\ldots,v_d)} z^v}
		{(1-z^{v_1})\cdot\ldots\cdot(1-z^{v_d})}.
	\label{eqn:rat-fun-cone}
\end{eqnarray}
If we specialize by substituting $z_i=q$ for each $i$, then we obtain $(1-q^\ell)^d$ in the denominator, since, by construction, all generators $v_i$ have coordinate sum $\ell$. This explains the appearance of binomial coefficients, since
\[
	\frac{1}{(1-q^\ell)^d} = \sum_{k=0}^\infty \binom{k+d-1}{d-1} q^{\ell\cdot k}
\]
which turns (\ref{eqn:counting-slices}) into
\begin{eqnarray}
	\sum_{k=0}^\infty \ehr_P(k)q^k = \frac{h^*_0 q^0 + \ldots + h^*_{d\cdot\ell-1}q^{d\cdot\ell-1}}{(1-q^\ell)^d}.
	\label{eqn:ehrhart-series}
\end{eqnarray}
In this way, many arithmetic calculations on the level $q$-series can be viewed as the projection of a geometric construction, via multivariate generating functions. The richer multivariate picture can be of use, for example, when converting arithmetic proofs into a bijective proofs, see \cite{BEKZ14}.

Intuitively, we can think of the generating functions $\phi_S$ as weighted indicator functions of sets of integer vectors. Starting with generating functions $\phi_P$ for polyhedra $P$ and taking linear combinations of these, we obtain an algebra $\PPP$ of polyhedral sets. However, working with rational function representations introduces an equivalence relation on this algebra. For example, we can expand $\frac{1}{1-z}$ either as $\sum_{i=0}^\infty z^i$, the indicator function of all non-negative integers, or as $-\sum_{i=1}^\infty z^{-i}$, minus the indicator function of all negative integers. This phenomenon generalizes to multivariate generating functions: To determine the formal expansion of a rational function uniquely, we have to fix a ``direction of expansion'' which can be given for example in terms of a suitable pointed cone. For details we refer the reader to, e.g., \cite{Monforte2013,Barvinok2008,Beck2009,Beck2007}. Important for our purposes is that to each rational function there corresponds an equivalence class of indicator functions and the simple example of the geometric series tells us what the equivalence relation is: Two elements in the algebra of polyhedral sets are equivalent if they are equal \emph{modulo lines}, i.e., modulo sets of the form $\mset{u+\lambda v}{\lambda\in\ZZ}$ for some $u,v\in\ZZ^n$. We say that a generating function $\phi$ is represented by some rational function expression $\rho$ if there exists a pointed cone $C$ such that the expansion of $\rho$ in the direction $C$ gives $\phi$; for this to be feasible we assume that the support of $\phi$ does not contain a line. Choosing a different direction $C'$ for the expansion of $\rho$ produces a generating function $\phi'$ that is equal to $\phi$ modulo lines.

\begin{figure}[t]
\includegraphics[angle=270,width=12.2cm]{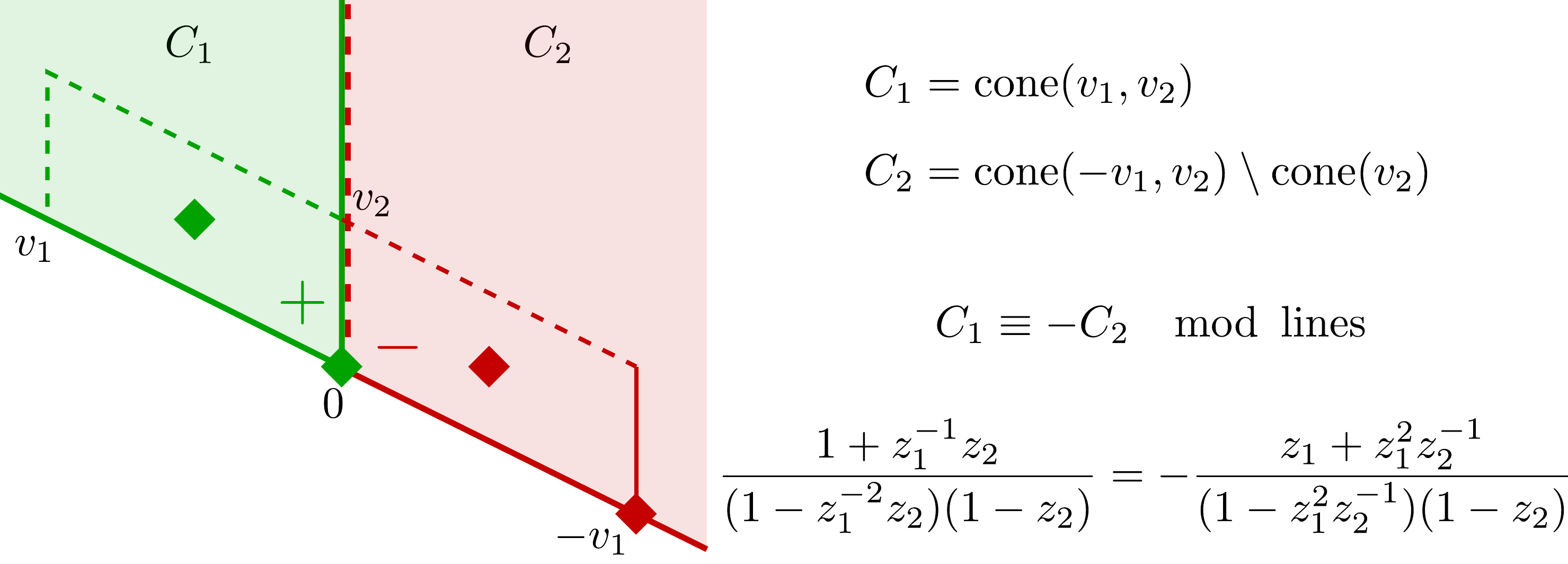}
\caption{Let $v_1 =(-2,1)$ and $v_2 =(0,1)$. Modulo lines, the closed cone $C_1$ generated by $v_1$ and $v_2$ is equal to the negative of the half-open cone $C_2$ generated by $-v_1$ and $v_2$. Integer points in corresponding fundamental parallelepipeds are shown as diamonds.}
\label{fig:flipping}
\end{figure}

Working with indicator functions of cones modulo lines does have its advantages. Most importantly, this allows us to ``flip'' cones by reversing the direction of some (or all) of their generators and opening some of their faces accordingly, as shown in Figure~\ref{fig:flipping}.

\begin{figure}[t]
\includegraphics[angle=270,width=12.2cm]{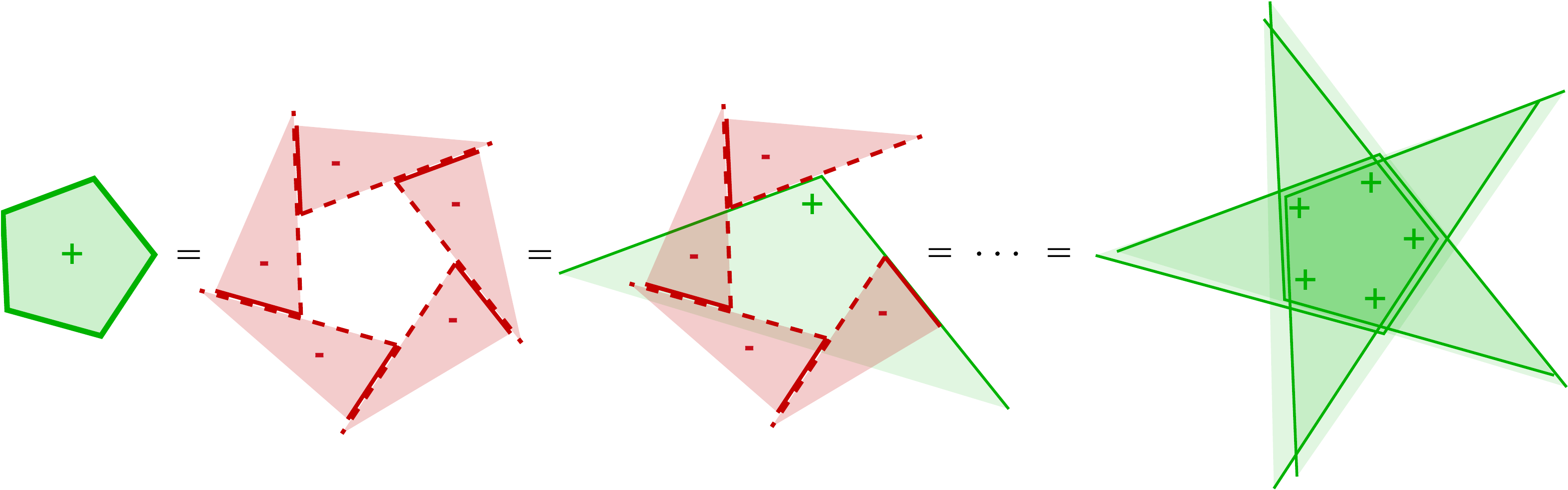}
\caption{Modulo lines, a polytope is equal to the sum of its vertex cones. In 2 dimensions, this is easy to see by iteratively flipping cones.}
\label{fig:brion}
\end{figure}

One beautiful application of this phenomenon is Brion's theorem, which allows us to represent $\phi_P$ for any line-free polyhedron $P$ in terms of rational function representations of cones, i.e., as a linear combination of expressions of the form (\ref{eqn:rat-fun-cone}). Brion's theorem is motivated in Figure~\ref{fig:brion}.

For a polyhedron $P$ we define the \emph{vertex cone} $\vcone(v,P)$ at a vertex $v$ of $P$ as the set
\[
	\vcone(v,P) = v + \cone_\RR(v_1,\ldots,v_N),
\]
where the $v_i$ are the directions of the edges incident to $v$, oriented away from $v$. We can easily represent each vertex cone by a rational function: For a simplicial cone $C$ we \emph{define} $\rho_C$ as the rational function expression given in (\ref{eqn:rat-fun-cone}).\footnote{Here it is important to note that (\ref{eqn:rat-fun-cone}) works also for cones with an apex $v\not=0$: All we have to do is take the fundamental parallelepiped $\Pi$ to be rooted at $v$ instead of the origin. This simply amounts to translating the fundamental parallelepiped as defined in Section~\ref{sec:cones} by $v$.} For a non-simplicial cone $C$ we define $\rho_C$ as a linear combination of such expressions, given via a triangulation of $C$. Then, the generating function of the set of integer points in $P$ is the sum of the rational function representations of the vertex cones.

\begin{theorem}[Brion \cite{Brion1988}]
Let $P$ be a polyhedron that does not contain any affine line. Then
\[
  \phi_P = \sum_{\text{$v$ vertex of $P$}} \rho_{\vcone(v,P)}(z).
\]
\end{theorem}

\begin{figure}[t]
\begin{center}
\includegraphics[angle=270,width=5.5cm]{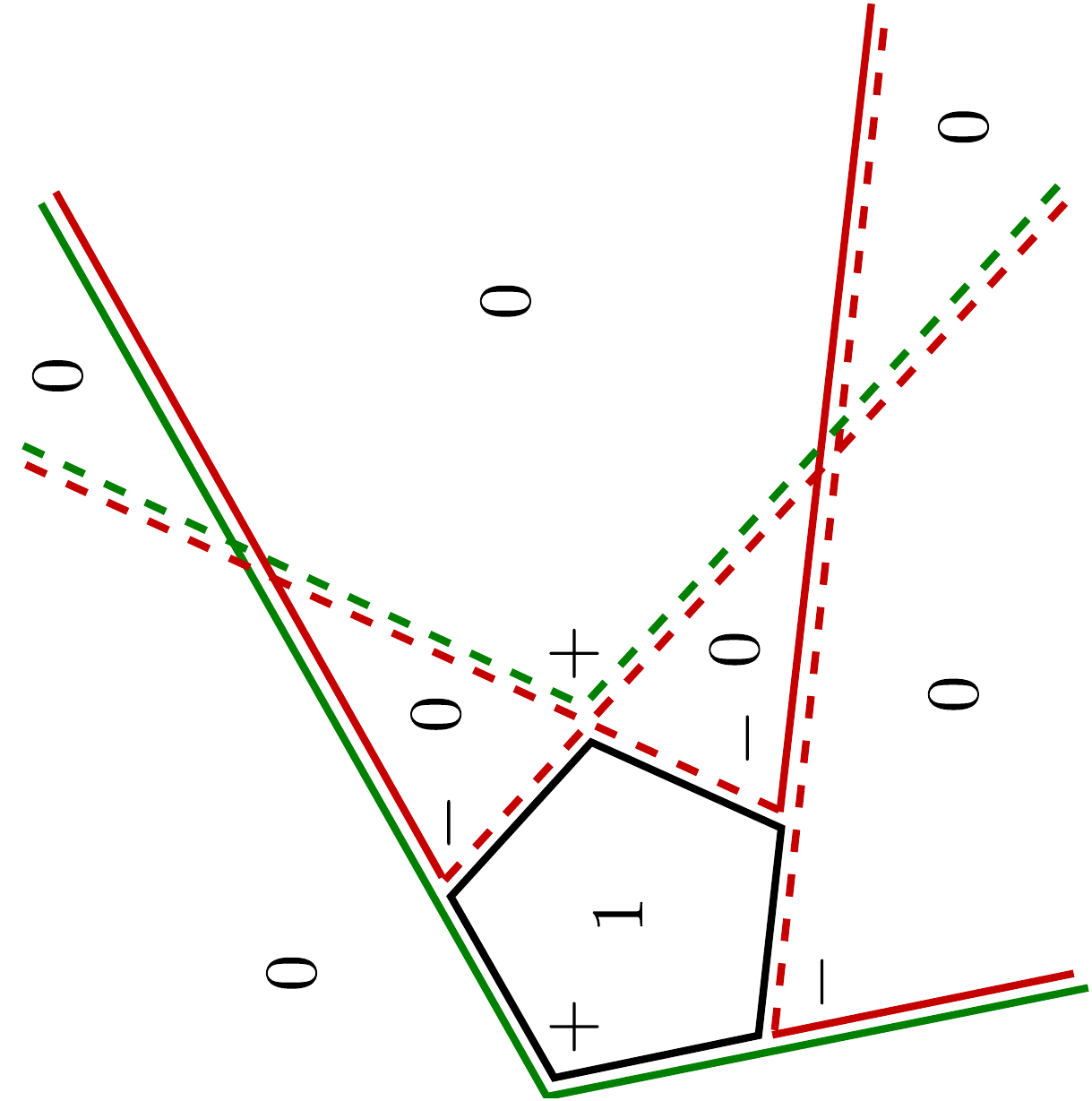}
\end{center}
\caption{The Lawrence-Varchenko decomposition of a pentagon. The sign next to the apex of each vertex cone $C$ specifies whether $C$ is to be added or subtracted. The number in each region is the net balance of how often points in the region are counted when the signed vertex cones are summed. All generators point of all vertex cones point to the right, which means all vertex cones are forward.}
\label{fig:lawrence}
\end{figure}

The theorem of Lawrence-Varchenko \cite{Lawrence1988,Varchenko1987} is the corresponding analogue for cases in which it is necessary to work with indicator functions directly, not with equivalence classes modulo lines. It expresses $\phi_P$ as an inclusion-exclusion of vertex cones which have been ``flipped forward'' so that their generators all point consistently in one direction of expansion as shown in Figure~\ref{fig:lawrence}.

\section{Coefficients of (Quasi-)Polynomials}
\label{sec:coefficients}

The geometric perspective provides a wide range of methods for establishing bounds on the coefficients of counting (quasi-)polynomials. In this section we will focus on polynomials for simplicity, but the results generalize to quasipolynomials.

The monomial basis is of course the classic choice for computing coefficients of polynomials. Geometrically, the elements of the monomial basis of the space of polynomials are the Ehrhart functions $\ehr_{[0,1)^i}(k)=k^i$ of half-open cubes $[0,1)^i$ of varying dimension. For us, it will be expedient to work with two different binomial bases instead, whose elements are the Ehrhart functions $\ehr_{\Delta^d_i}(k) = \binom{k+d-i}{d}$ of unimodular\footnote{A simplex $\Delta$ with integer vertices is \emph{unimodular} if the fundamental parallelepiped of $\cone(\Delta\times\{1\})$ contains only a single integer vector: the origin. Equivalently $\ZZ^n\cap\cone_\RR(\Delta\times\{1\}) = \cone_\ZZ(\Delta\times\{1\})$.} $d$-dimensional \emph{half-open simplices} $\Delta^d_i$ with $i$ open facets. Up to lattice equivalence, such a $\Delta^d_i$ has the form
\[
  \Delta^d_i = \mset{x\in\RR^{d+1}}{x_1 > 0,\ldots,x_i > 0, x_{i+1} \geq 0, \ldots, x_{d+1} \geq 0, \sum_j x_j = 1}.
\]

These unimodular half-open simplices $\Delta^d_i$ form the basic building block of Ehrhart theory. They offer two different ways in which we can use them to construct a basis of the space of polynomials. The first basis, which defines the $h^*$-coefficients, fixes the dimension $d$ of the simplices and varies the number $i$ of open facets. In contrast, the second basis, which defines the $f^*$-coefficients, uses only open simplices with $i=d+1$, but varies their dimension $d$.

Formally, the $h^*$-vector $(h^*_0,\ldots,h^*_d)$ and the $f^*$-vector $(f^*_0,\ldots,f^*_d)$ of a polynomial $p(k)$ of degree at most $d$ are defined by
\begin{eqnarray*}
p(k) 
& = & 
h^*_0 \binom{k+d}{d} +  h^*_1 \binom{k+d-1}{d} + \ldots + h^*_d \binom{k}{d} \\
& = &
f^*_0 \binom{k-1}{0} +  f^*_1 \binom{k-1}{1} + \ldots + f^*_d \binom{k-1}{d}.
\end{eqnarray*} 

\begin{figure}[t]
\includegraphics[angle=270,width=12.2cm]{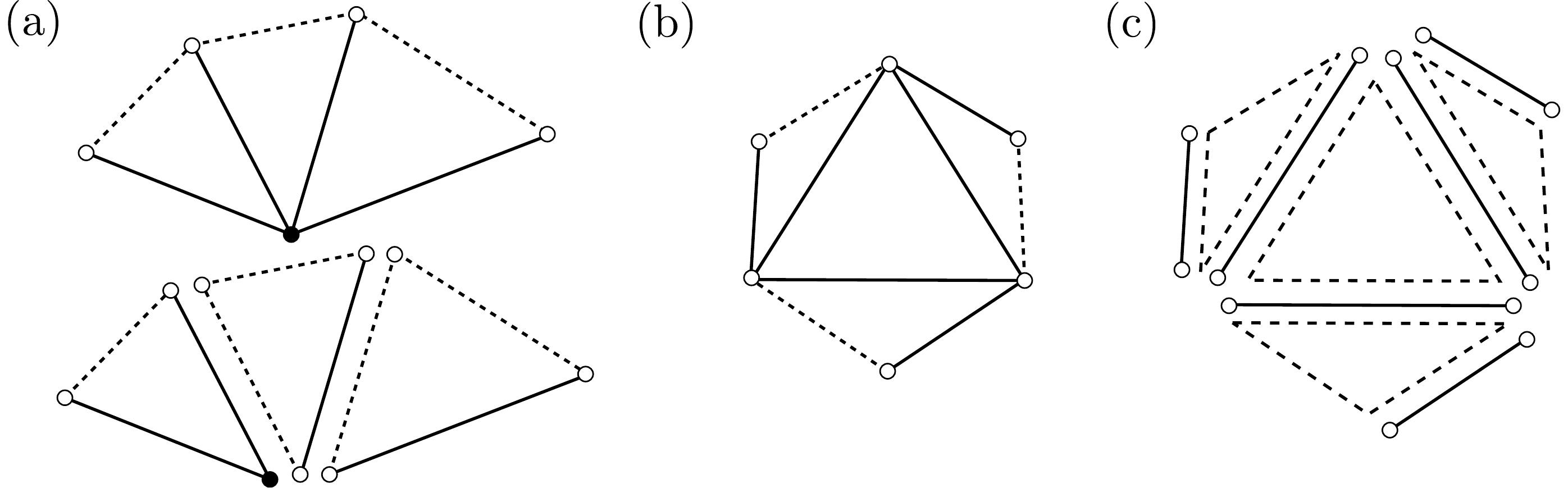}
\caption{(a) A half-open 2-dimensional partial polytopal complex and a partition into half-open 2-dimensional simplices. (b) A half-open partial polytopal complex $X$ that is not partitionable. A partition of this complex with half-open 2-dimensional simplices would have to contain a 2-dimensional simplex with at least two edges and, consequently, at least one vertex. However, $X$ does not contain any of its vertices. (c) A decomposition of $X$ into open simplices of various dimension.}
\label{fig:partitionability}
\end{figure}

Let us begin by taking a closer look at the $h^*$-coefficients. As we have seen in (\ref{eqn:counting-slices}) and (\ref{eqn:ehrhart-series}), the $h^*$-vector of the Ehrhart quasipolynomial of a simplex $\Delta$ counts lattice points at different heights in the fundamental parallelepiped of $\cone(\Delta\times\{1\})$, which immediately implies $h^*_i\geq 0$. This observation extends to half-open simplices where some facets have been removed. It follows that if a geometric model $X$ can be \emph{partitioned} into half-open simplices that are all of full dimension, as shown in Figure~\ref{fig:partitionability}(a), it follows that $\ehr_X$ has non-negative $h^*$-vector as well. As it turns out, all (closed convex) polytopes have such a \emph{partitionable} triangulation, which proves non-negativity of the $h^*$-vector for all polytopes.

\begin{theorem}[Stanley \cite{Stanley1980}]
If $P$ is an integral polytope, then $\ehr_P$ has a non-negative $h^*$-vector.
\end{theorem}

However, as the examples of the chromatic polynomial and the flow polynomial from Sections~\ref{sec:modeling} and \ref{sec:reciprocity} show, the geometric models $X$ that appear in combinatorial applications of Ehrhart theory are not simply polytopes: Often they are non-convex, disconnected, half-open or have non-trivial topology. This can lead to geometric models $X$ that are not partitionable and, consequently, to counting polynomials with negative $h^*$-coefficients. 

Figure~\ref{fig:partitionability}(b) gives an example of a half-open partial polytopal complex that is not partitionable: A partition of the complex in Figure~\ref{fig:partitionability}(b), for example, would require 4 half-open simplices of dimension 2 that have, in total, 6 closed edges but contain none of the vertices of the complex, which is impossible. Here it is important to recall that, because we are working with the $h^*$ basis, all half-open simplices participating in a partition are required to have the same dimension (in this example, dimension 2).

Such phenomena appear in practice. One prominent example of natural counting polynomials with negative entries in their $h^*$-vector are chromatic polynomials of hypergraphs. In this case, it is the non-trivial topology of the geometric models that gives rise to non-partitionability: It is easy to construct hypergraphs whose coloring complexes consists of, say, 2-dimensional spheres that intersect in 0-dimensional subspheres; such complexes are not partitionable and can produce negative $h^*$-coefficients \cite{BDK2012}. 

As we have seen in Section~\ref{sec:ehrhart}, partial polytopal complexes are the right notion to describe combinatorial models in Ehrhart theory. While partial polytopal complexes are not always partitionable, they can always be written as a disjoint union of relatively open simplices of various dimension. The partial polytopal complex in Figure~\ref{fig:partitionability}(b) can, for example, be written as a disjoin union of open simplices of dimension 1 and 2 as shown in Figure~\ref{fig:partitionability}(c). This motivates the use of the $f^*$-basis. As it turns out, the $f^*$-vector of an open simplex $\Delta$ has a counting interpretation similar to (\ref{eqn:counting-slices}), even though its construction is more subtle \cite{Breuer2012}. It follows that all partial polytopal complexes with integer vertices have a non-negative $f^*$-vector. Moreover, this property characterizes Ehrhart polynomials of partial polytopal complexes.

\begin{theorem}[Breuer \cite{Breuer2012}]
\label{thm:f-star}
If $X$ is an integral partial polytopal complex, then $\ehr_X$ has a non-negative $f^*$-vector. 

Conversely, if $p(k)$ is a polynomial with non-negative $f^*$-vector, then there exists an integral partial polytopal complex $X$ such that $\ehr_X(k)=p(k)$.
\end{theorem}

While Theorem~\ref{thm:f-star} characterizes Ehrhart polynomials of the kind of geometric objects that appear in many combinatorial applications, the question remains how to characterize Ehrhart polynomials of convex polytopes. This challenge is vastly more difficult, and, even though many constraints on the $h^*$-vectors of convex polytopes have been proven, is still wide-open even in dimension 3. At least in dimension 2, a complete characterization of the coefficients of Ehrhart polytopes is available. See \cite{Beck2005,Haase2009,Henk2009,Stapledon2009} for more information.

Still, there are a wealth of tools available for proving sharper bounds on the coefficients of counting polynomials $\ehr_X$, by exploiting the particular geometric structure of the partial polytopal complex $X$, even if $X$ is not convex. One of the most powerful techniques available is the use of convex ear decompositions. A \emph{convex ear decomposition} is a decomposition of a simplicial complex $X$ into ``ears'' $E_0,\ldots,E_N$ such that $E_0$ is the boundary complex of a simplicial polytope, the remaining $E_i$ are balls that are subcomplexes of the boundary complex of some simplicial polytope, and $E_i$ is attached to $\bigcup_{j<i}E_j$ along its entire boundary (and not just along some facets), i.e., $E_i \cap \bigcup_{j<i}E_j = \partial E_i$. For example, the complex in Figure~\ref{fig:chromatic}, consisting of the boundary of the cube and the two hyperplanes, has a convex ear decomposition: Start with the boundary of the cube as triangulated by the braid arrangement, glue in the triangulated square lying on one of the hyperplanes and then glue in the two triangles on the second hyperplane one after the other. If all simplices in this complex are unimodular (as in many combinatorial applications), this leads to the following bounds, which have been successfully applied to the chromatic polynomial by Hersh and Swartz \cite{Hersh2007} and to the integral and modular flow and tension polynomials by Breuer and Dall \cite{Breuer2011}.

\begin{theorem}[Chari \cite{Chari1997}, Swartz \cite{Swartz2006}]
\label{thm:convex-ear-decomposition}
If $X$ is a simplicial complex in which all simplices are unimodular and $X$ has a convex ear decomposition then the $h^*$-vector of $\ehr_X(k)$ satisfies
\begin{enumerate}[a)]
	\item $h^*_0 \leq h^*_1 \leq \cdots \leq h^*_{\floor{d/2}}$,
	\item $h^*_i \leq h^*_{d-i}$ for $i\leq d/2$, and
	\item $(h^*_0,h^*_1-h^*_0,\ldots,h^*_{\ceil{d/2}} - h^*_{\ceil{d/2}-1})$ is an $M$-vector\footnote{$M$-vectors are defined as in Macaulay's theorem, see for example \cite[Chapter~8]{Ziegler1995}.}.
\end{enumerate}
\end{theorem}

\section{Quasisymmetric Functions}
\label{sec:quasisymmetric}

Polyhedral models are useful for the study of combinatorial objects beyond counting polynomials as well. For example, the simple construction from Section~\ref{sec:modeling} of intersecting the cube with a subarrangement of the braid arrangement can serve as a lens into the world of quasisymmetric functions \cite{Breuer2014}.

A \emph{quasisymmetric function} is a formal power series $Q$ of bounded degree in countably many variables $x_1,x_2,\ldots$ such that the coefficients of $Q$ are shift invariant, i.e., for every $(\alpha_1,\ldots,\alpha_m)$ the coefficients of the monomials $x_{i_1}^{\alpha_1}x_{i_2}^{\alpha_2}\cdots x_{i_m}^{\alpha_m}$ for any $i_1<i_2<\ldots<i_m$ are equal \cite{Stanley2001}. Note that a quasisymmetric function can have bounded degree without being a polynomial since we have infinitely many variables at our disposal.

To approach these from a geometric perspective, it is instructive to start with quasisymmetric functions in non-commuting variables or \emph{nc-quasisymmetric functions} for short \cite{Bergeron2009}. Here the variables $x_i$ do not commute multiplicatively and the constraint is that two monomials $x_{i_1}\cdots x_{i_d}$ and $x_{j_1}\cdots x_{j_d}$ have the same coefficient if the tuples $i=(i_1,\ldots,i_d)$ and $j=(j_1,\ldots,j_d)$ induce the same ordered set partition $\Delta(i)=\Delta(j)$. Here $\Delta(i)=(\Delta_1,\ldots,\Delta_m)$ is an ordered partition of the index set $\{1,\ldots,d\}$ such that $i|_{\Delta_l}$ is constant and $i|_{\Delta_l}<i|_{\Delta_{l+1}}$ for all $l$, e.g., $\Delta(3,2,2,3,1) = (\{5\},\{2,3\},\{1,4\})=:5|23|14$.

\begin{figure}[t]
\includegraphics[angle=270,width=12.2cm]{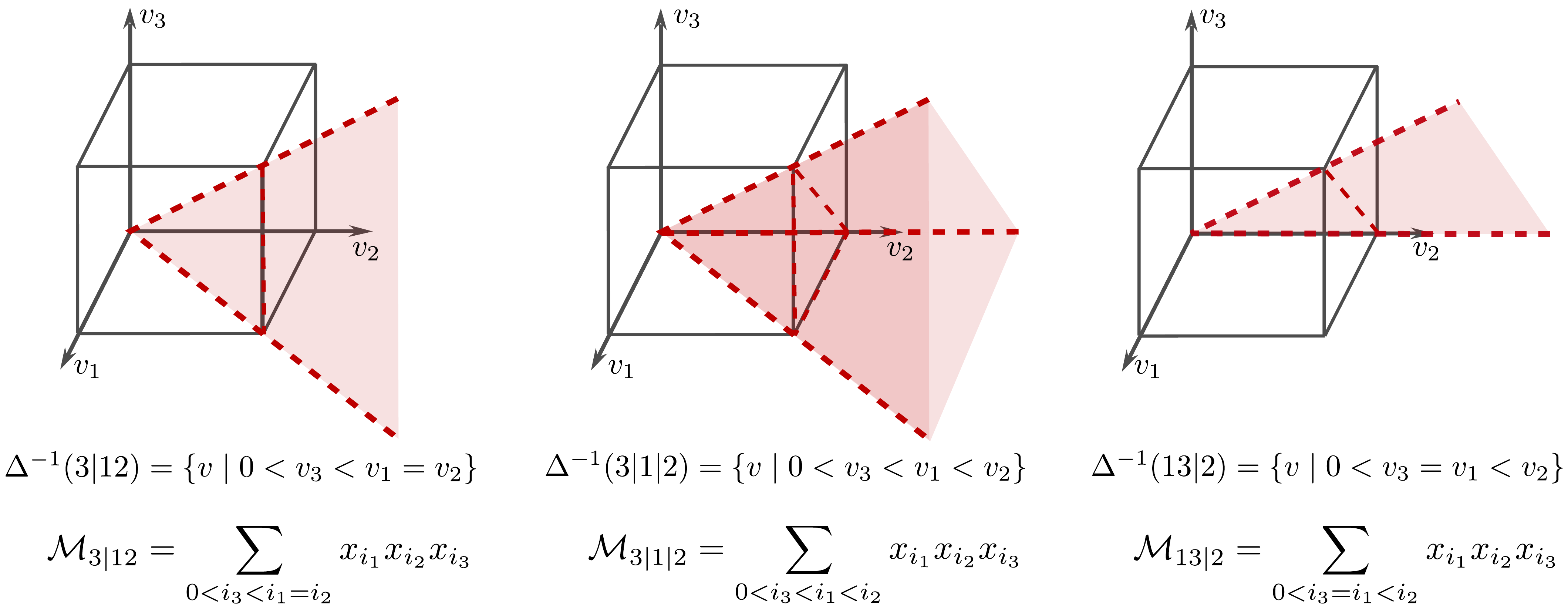}
\caption{The ``shift-invariant regions'' of integer points that are mapped to the same ordered set partition by $\Delta$ are the simplicial cones in the braid arrangement and correspond to the monomial quasisymmetric function.}
\label{fig:quasisymmetric}
\end{figure}

To visualize what is going on here, we need a new way of associating integer vectors with monomials. Classically, we identify monomials in commuting variables with their exponent vector. Here, we identify monomials in non-commuting variables with their vector of indices, i.e., we identify $x_{v_1}\cdots x_{v_d}$ with $(v_1,\ldots,v_d)\in\ZZ^d_{\geq 1}$. This allows us to picture the map $\Delta$: If $\phi$ is an ordered set partition of $\{1,\ldots,d\}$, then $\Delta^{-1}(\phi)$ is precisely the set of integer vectors contained in a simplicial cone of the partial polyhedral complex obtained by triangulating the positive orthant by the braid arrangement, as shown in Figure~\ref{fig:quasisymmetric}. In other words, the \emph{monomial nc-quasisymmetric functions}
\[
	\MMM_\phi = \sum_{v\in\ZZ^d, \Delta(v)=\phi} x_{v_1}\cdots x_{v_d}
\]
form a basis of the space of nc-quasisymmetric functions and these are nothing but cones in the braid arrangement.

Any nc-quasisymmetric function can be turned into a quasisymmetric function simply by allowing variables to commute. This can be modeled geometrically by taking an integer vector and permuting its entries so that they are in weakly increasing order, i.e., an element of the half-open simplicial cone $C:=\mset{v}{0<v_1\leq\ldots\leq v_d}$. This maps $\MMM_{(\phi_1,\ldots,\phi_m)}$ to the \emph{monomial quasisymmetric function} $M_{(|\phi_1|,\ldots,|\phi_m|)}$ where
\[
  M_{(\alpha_1,\ldots,\alpha_m)} = \sum_{1\leq i_1<\ldots<i_m} x_{i_1}^{\alpha_1}\cdot\ldots\cdot x_{i_m}^{\alpha_m}.
\]
The monomial quasisymmetric functions form a basis of the space of quasisymmetric functions. Thus every quasisymmetric function can be visualized as assigning a weight to every face of the cone $C$. The support of a quasisymmetric function is thus a partial polyhedral subcomplex $X$ of the face lattice of $C$.

Going one step further it is possible to obtain a polynomial $p$ from a quasisymmetric function $Q$ by substituting 1 into the first $k$ variables and $0$ into all other variables, i.e., $Q(1^k)=p(k)$. Geometrically, this substitution eliminates all integer points that contain an entry larger than $k$. This corresponds to intersecting the complex $X$ of cones with the cube $(0,k]^d$, turning $X$ into a simplicial complex $X\cap(0,k]^d$ and $p$ into the Ehrhart function $\ehr_{X\cap(0,1]^d}(k) = Q(1^k)$.\footnote{This works best if $Q$ is the specialization of an nc-quasisymmetric function with 0-1 coefficients. Otherwise, this would require a linear combination of Ehrhart functions.} These observations provide a direct translation between Ehrhart functions constructed using the braid arrangement and quasisymmetric functions.

This connection provides fertile ground for future exploration. On the one hand, the geometric approach offers a very flexible framework for defining quasisymmetric functions. Scheduling problems alone capture a wide range of known quasisymmetric functions, such as the chromatic symmetric function, the matroid invariant of Billera-Jia-Reiner, or Ehrenborg's quasisymmetric function for posets, as well as new ones, such as the Bergman and arboricity quasisymmetric functions \cite{Breuer2014}. On the other hand, many methods for the analysis of Ehrhart polynomials carry over to the quasisymmetric function world. For example, the specialization $Q(1^k)$ collects the coefficients of $Q$ in the fundamental basis in the $h^*$-vector of the associated Ehrhart-polynomial -- and similarly for the monomial basis and the $f^*$-vector. In particular, if $X$ is a partial subcomplex of the braid arrangement and $N, Q$ and $\ehr_{X\cap(0,1]^d}$ are the associated nc-quasisymmetric, quasisymmetric and Ehrhart functions, then partitionability of $X$ implies non-negativity of the coefficients in the fundamental basis of $N$ and $Q$ and non-negativity of the $h^*$-vector of $\ehr_{X\cap(0,1]^d}$. If $X$ is given by a scheduling problem, partitionability can be guaranteed if the boolean expression defining the scheduling problem takes the form of a certain kind of decision tree \cite{Breuer2014}.

\section{Algorithms for Counting Integer Points in Polyhedra}
\label{sec:computing}

There are many different computational problems associated with polyhedra. The problem of deciding whether there exists a rational vector $v\in\QQ^n$ satisfying a linear system of inequalities\footnote{Solving a linear system of \emph{inequalities} over $\ZZ$ (or, equivalently, solving a linear system of equations over $\NN$) is NP-hard. However, solving a linear system of \emph{equations} over $\ZZ$ is polynomial-time solvable, for example using the Smith normal form, see below.} is polynomial time computable, but when we look for an integer vector $v\in\ZZ^n$ instead, the problem becomes NP-hard \cite{Schrijver1986}. However, if the dimension of the polyhedron, i.e., the number of variables of the system, is fixed a priori, then there is a polynomial time algorithm for finding an integer solution as Lenstra was able to show in 1983 \cite{Lenstra1983}. While the problem of counting integer solutions is \#P-hard as well, the question remained open whether it becomes polynomial time computable if the dimension is fixed. The first algorithm with a polynomial running time in fixed dimension was described by Barvinok in 1994 \cite{Barvinok1994} and it took ten more years until such an algorithm was first implemented by De Loera et al. in 2004 \cite{DeLoera2004}. 

In this section we give an overview over the algorithmic methods for computing the number of integer points in a polyhedron $P$, and the related problems of computing the Ehrhart polynomial $\ehr_P$ and a rational function expression of the multivariate generating function $\phi_P$ of all integer points in $P$. Independently of whether the goal is to compute $\ehr_P$ by first passing from $P$ to $\cone(P\times\{1\})$ or whether the goal is to compute $\phi_P$ and $\ZZ^n\cap P$ directly by using Brion's theorem, the methods employed are similar and consist of three basic steps. First, the polyhedron $P$ is decomposed into simplicial cones. Second, a rational function representation of the integer points in these simplicial cones is computed. We will focus on this step in our exposition since it is crucial with regard to runtime complexity. Third, the obtained rational function expression needs to be specialized if the number of integer points or the Ehrhart (quasi-)polynomial is desired.

To decompose a polyhedron $P$ into simplicial cones, we start by appealing to Brion's theorem and represent $P$ as the sum of its vertex cones, modulo lines.\footnote{We can also use the theorem of Lawrence-Varchenko to obtain an exact signed decomposition, without working modulo lines.} To achieve this, we need to compute the vertices and edge directions of $P$. Next, the resulting cones need to be triangulated to make them simplicial. There are sophisticated algorithms available for both tasks \cite{Fukuda1994,Fukuda1996,Pfeifle2003,DeLoera2010}. It is also possible to compute a decomposition of $P$ into simplicial cones directly, without computing vertices or triangulating, using the Polyhedral Omega algorithm \cite{BZ14}. Polyhedral Omega is based on simple explicit rules for manipulating simplicial cones formally and is motivated by the symbolic computation framework of partition analysis \cite{Andrews2001}.

\begin{figure}[t]
\begin{center}
\includegraphics[angle=270,width=12.2cm]{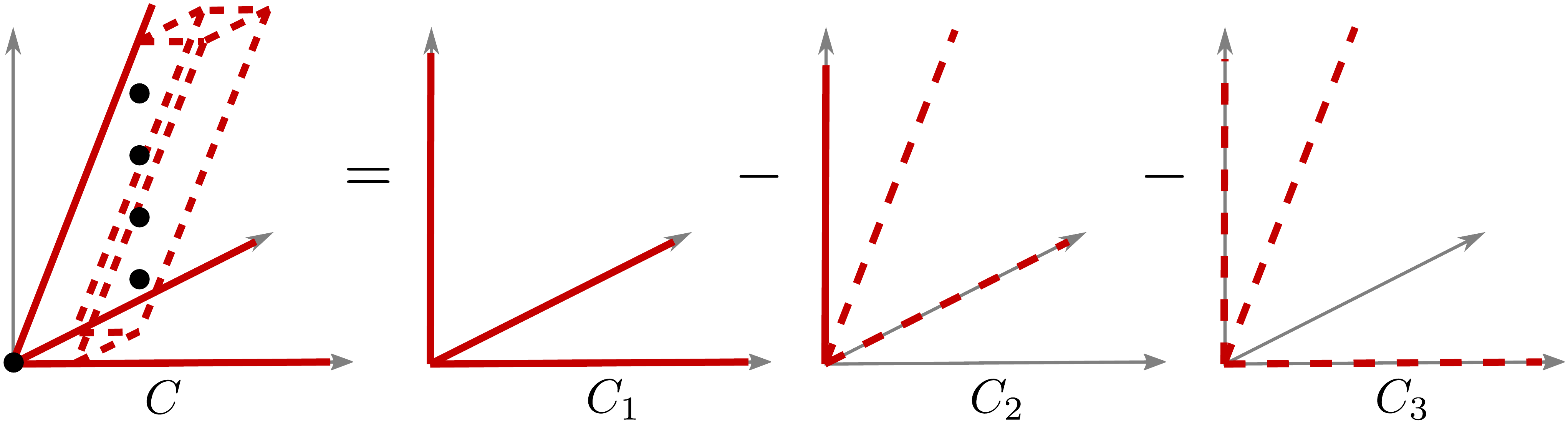}
\end{center}
\caption{$C=\cone_\RR((1,0,0),(0,1,0),(1,1,a))$ has $a$ integer points in its fundamental parallelepiped. This number of integer points is therefore exponential in the encoding size of $C$ which is in $\OOO(\log(a))$. However $C$ can be written as a signed sum $C=C_1-C_2-C_3$ of unimodular cones. Here the facet of $C_2$ generated by $(0,1,0)$ and $(1,1,a)$ is open and the two facets of $C_3$ generated by $(1,1,a)$ and one of the other two generators are open.}
\label{fig:barvinok}
\end{figure}

In the second step, we use the ideas developed in Sections~\ref{sec:cones} and \ref{sec:brion} to represent the generating function $\phi_C$ of integer points in a simplicial cone $C=\cone_\RR(v_1,\ldots,v_d)\subset\ZZ^d$ as a rational function. Let $V$ denote the matrix with the generators $v_i$ of $C$ as columns. The straightforward approach is to use (\ref{eqn:rat-fun-cone}) and obtain a rational function expression by simply enumerating all integer points in the fundamental parallelepiped $\Pi(V)$ of $C$. This is both simple and efficient if the index $\ZZ^d\cap \Pi=|\det(V)|$ of $C$ is sufficiently small. However, in the worst case, the index may be exponential in encoding size of $V$, as Figure~\ref{fig:barvinok} shows. Thus it is not clear a priori that there exists a rational function expression for $\phi_C$ whose encoding size is polynomial in the encoding size of the input. Barvinok's key achievement was to find such a representation. 

\begin{figure}[t]
\begin{center}
\includegraphics[angle=270,width=12cm]{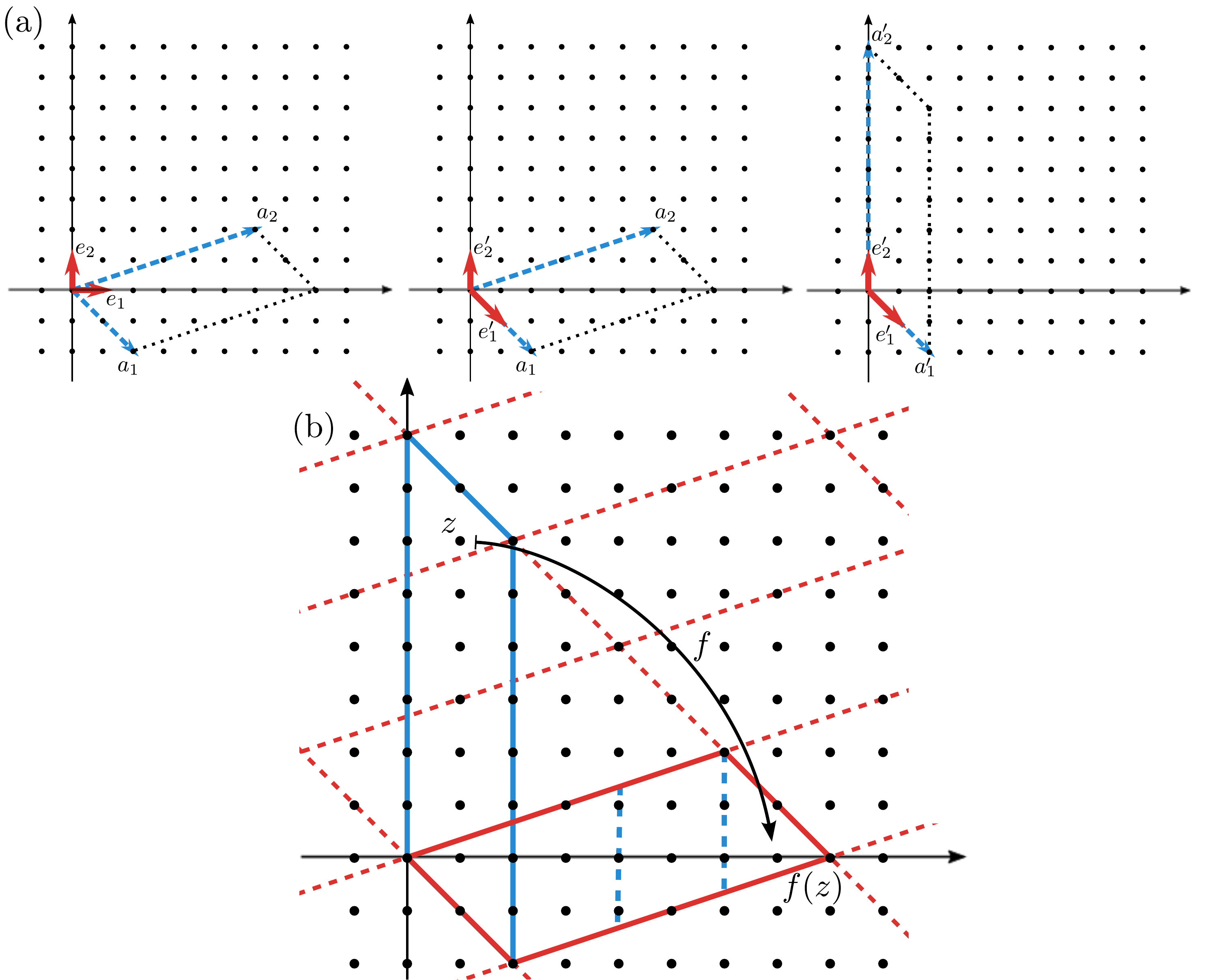}
\end{center}
\caption{In order to compute all integer points in the fundamental parallelepiped $\Pi(a_1,a_2)$, shown in the left panel of (a), we proceed as follows. (a) Using the Smith normal form, we first perform a change of basis on the integer lattice $\ZZ^2$ and then a change of basis on the sublattice generated by $a_1,a_2$, so that the bases align. (b) Listing all the integer points in the aligned fundamental parallelepiped $\Pi(a_1',a_2')$ is easy. We transform these into integer points in $\Pi(a_1,a_2)$ by modular arithmetic (taking fractional parts of coordinates wrt.~the original basis $a_1,a_2$).}
\label{fig:snf}
\end{figure}

Before we come to Barvinok's short rational function representation, however, it is instructive to take a closer look at how to enumerate the integer points in $\Pi$ explicitly. There are several well-known approaches to this problem \cite{BZ14,BIS2012,Koppe2008} which are all closely related. We will work with the Smith normal form of the matrix $V$, which can be computed in polynomial time \cite{Schrijver1986}. The Smith normal form of $V$ is a representation $V=USW$ where $U,S,W$ are integer matrices, $U,W$ have determinant $\pm1$ and $S$ is a diagonal matrix whose diagonal entries $s_1,\ldots,s_d$ satisfy $s_i|s_{i+1}$. This can be interpreted as shown in Figure~\ref{fig:snf}. The columns of $V$ form a basis of a sublattice $J$ of the integer lattice $\ZZ^d$, and $V$ gives the coordinates of this basis with respect to the standard basis of $\ZZ^d$. The matrices $U$ and $W$ represent changes of basis on \emph{both} lattices such that the new bases $B_J$ of $J$ and $B_{\ZZ^d}$ of $\ZZ^d$ line up. Since the elements of $B_J$ are multiples of the elements of $B_{\ZZ^d}$, the integer points $x_i$ in the fundamental parallelepiped of $B_J$ are easy to enumerate. By computing the coordinates of the $x_i$ wrt.~the original basis $V$ of $J$ and taking fractional parts, we translate the $x_i$ into the fundamental parallelepiped $\Pi(V)$ and we are guaranteed that we get every point in $\ZZ^d\cap\Pi(V)$ exactly once. This process is summarized in the formula
\[
  \phi_{\cone_\RR(V)(z)} = \frac{\sum_{k_1=0}^{s_1-1}\cdots\sum_{k_1=0}^{s_1-1} z^{\frac{1}{s_d} V( W^{-1} (s_1'k_1,\ldots,s_d'k_d)^\top \mod s_d ) } }{(1-z^{v_1})\cdot\ldots\cdot(1-z^{v_d})}
\]
where $s_i'=\frac{s_d}{s_i}$. This particular expression is taken from \cite{BZ14}.

Now we come to Barvinok's central idea. Consider the cone $C$ generated by $(1,0,0)$, $(0,1,0)$ and $(1,1,a)$ for $0<a\in\ZZ$. Its fundamental parallelepiped contains $a$ integer points, as shown in Figure~\ref{fig:barvinok}, which is exponential in the encoding size $\OOO(\log(a))$ of $C$. Moreover, there is no way to write $C$ as a union of $\OOO(\log(a))$ unimodular cones of index 1. Using inclusion-exclusion, however, $C$ can be written as the positive orthant $C_1$ minus the cone $C_2$ generated by $(0,0,1)$, $(0,1,0)$, $(1,1,a)$ and the cone $C_3$ generated by $(1,0,0)$, $(0,0,1)$, $(1,1,a)$ which all have index 1. This generalizes. Let $C$ denote a simplicial cone in fixed dimension $d$ and let $I$ denote its index. Using the LLL algorithm it is possible to find an integer vector $u$ such that $C=\cone_\RR(v_1,\ldots,v_d)$ can be written as a signed combination of the cones $C_1=\cone_\RR(u,v_2\ldots,v_d)$, $C_2=\cone_\RR(v_1,u,\ldots,v_d)$, $\ldots$, $C_d=\cone_\RR(v_1,\ldots,v_{d-1},u)$, where some facets of the $C_i$ have to be opened according to a few explicit combinatorial rules \cite{Koppe2008}. The key property of this construction is that indices of the cones $C_i$ decrease quickly. Applying this decomposition recursively, the indices of the cones will eventually reach 1, i.e., the cones will become unimodular. At each node of the recursion tree one cone is split into $d$-cones, however, the depth of the tree is at most doubly logarithmic in $I$. Thus the total number of cones obtained is polynomial in the encoding length of $C$. The result is the following fundamental theorem.

\begin{theorem}[Barvinok \cite{Barvinok1994}]
Let $C\subset\ZZ^d$ be a $d$-dimensional simplicial cone with integer generators. Then there exists signs $\epsilon_i$ and vectors $a_i,b_{i,j}$ such that
\begin{eqnarray}
  \phi_C(z) = \sum_{i=1}^N \epsilon_i\frac{z^{a_i}}{(1-z^{b_{i,1}})\cdot\ldots\cdot(1-z^{b_{i,d}})}
\end{eqnarray}
and for fixed $d$ the number of summands $N$ is bounded by a polynomial in the encoding length of $C$.
\end{theorem}

The third step is to specialize the representation of $\phi_P$ in terms of multivariate rational functions we have obtained thus far, in order to get the Ehrhart polynomial $\ehr_P$ or the number $\#\ZZ^n\cap P$. This specialization is non-trivial, especially if Barvinok decompositions are used, since typically the desired specialization is a pole of the rational function representation. However, using an exponential substitution and limit arguments it is possible to compute this specialization in polynomial time.

The toolbox of algorithms we have described here has many more applications and extensions. For example, it is possible to extend these methods to handle multivariate Ehrhart polynomials \cite{Verdoolaege2007}, to compute intersections $\phi_{P\cap Q}$ given $\phi_P$ and $\phi_Q$ \cite{Barvinok2003}, to compute Pareto optima in multi-criteria optimization over integer points in polyhedra \cite{DeLoera2009}, to integrate and sum polynomials over polyhedra \cite{Baldoni2011} and to convert between rational function representations and piecewise quasipolynomial representations of counting functions \cite{Verdoolaege2008} -- all in polynomial time if the dimension is fixed. As starting points for further reading we recommend the textbooks \cite{Barvinok2008,DeLoera2012}.

\section{Lattice Point Sets and the Euclidean Algorithm}
\label{sec:euclid}

After these very general considerations, we end this exposition on a playful note by taking a closer look at integer point geometry in dimension 2 and discussing several different ways in which the Euclidean algorithm makes an appearance.

\begin{figure}[t]
\begin{center}
\includegraphics[angle=270,width=12.2cm]{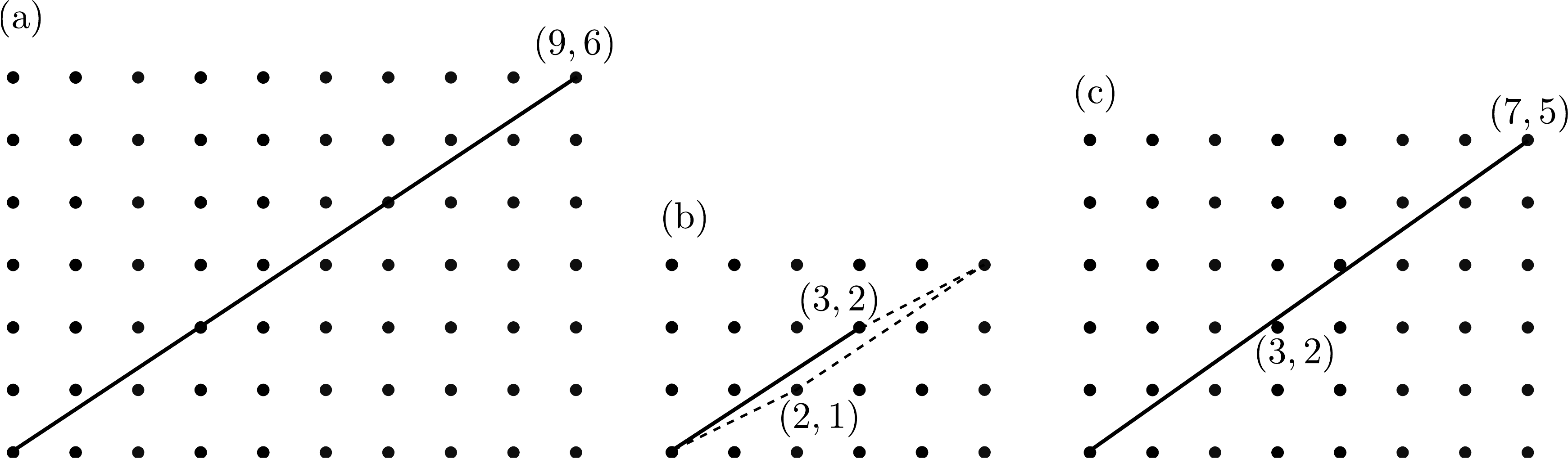}
\end{center}
\caption{(a) $\gcd(9,6)=3.$ (b) $-1\cdot 3 + 2 \cdot 2 = 1$ since $(2,1)$ is closest to the line through $(3,2)$ which means $\ZZ^2\cap\Pi((3,2),(2,1))$ contains no integer point except the origin. (c) $\gcd(7,5)=1$ and $-2\cdot 7 + 3\cdot 5 = 1$.}
\label{fig:gcd-plane}
\end{figure}

The integer lattice in the plane is a great stage for visualizing the greatest common divisor, as Figure~\ref{fig:gcd-plane} shows. For two integers $a,b\in\ZZ$, the line segment in the plane from the origin to the point $(a,b)$ contains precisely $\gcd(a,b)+1$ integer points. Let $(p,q)$ denote the coordinates of a lattice point closest to but not on the line $L$ through $(0,0)$ and $(a,b)$. By construction, the fundamental parallelepiped spanned by $(a,b)$ and $(p,q)$ contains precisely $\gcd(a,b)$ lattice points on the line and no lattice points off the line. 
\[
\gcd(a,b)=\Pi((a,b),(p,q))=\det\mat{a&p\\b&q} = ap - bq.
\]
Thus the coordinates of the closest points give precisely the coefficients produced by the extended Euclidean algorithm.

\begin{figure}[t]
\begin{center}
\includegraphics[angle=270,width=12.2cm]{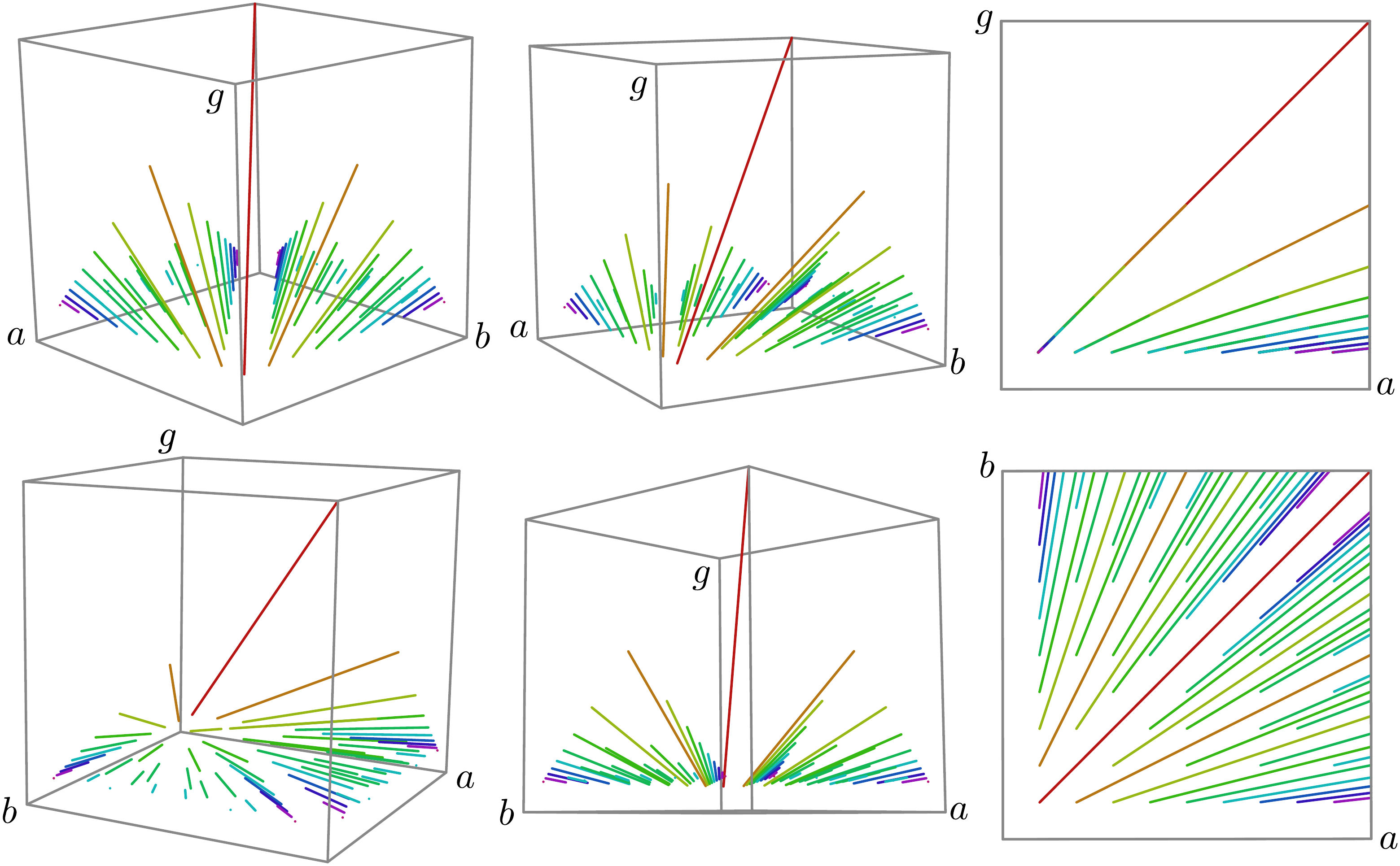}
\end{center}
\caption{The graph of the $\gcd$ $g=\gcd(a,b)$ as described in the text. Shown are the rays from the origin through $(a,b)$ starting at $(a,b)$ for all $a,b$ with $\gcd(a,b)=1$. The color of a ray is given by its depth in the recursion tree of the Euclidean algorithm. The plots in the rightmost column show parallel projections of the graph onto the $(g,a)$ and $(a,b)$ planes, respectively.}
\label{fig:gcd-graph}
\end{figure}

From the above observation it immediately follows that the value of the GCD increases linearly along any such line $L$. If $(a,b)$ is the integer point closest to the origin on such a line $L$, then $\gcd(a,b)=1$. The next values of the GCD on $L$ are thus $\gcd(2a,2b)=2$, $\gcd(3a,3b)=3$. The graph of the function $\gcd:\ZZ^2_{>0}\rar \ZZ_{>0}$ is thus contained in a countable collection of rays from the origin through all points $(a,b)$ with $\gcd(a,b)=1$. This ``graph'' of the GCD is shown in Figure~\ref{fig:gcd-graph}.

A closer look at the graph in Figure~\ref{fig:gcd-graph} immediately reveals a recursive tree-like structure. It turns out that this tree corresponds precisely to the recursive operation of the Euclidean algorithm. The Euclidean algorithm as described by Euclid moves from $(a,b)$ to $(a-b,a)$ if $a>b$, it moves from $(a,b)$ to $(a,b-a)$ if $a<b$ and it terminates if $a=b$. The perceptive reader will note that this immediately gives a way to enumerate all positive rational numbers as nodes of an infinite binary tree \cite{Calkin2000}. However, tracing out these paths of the Euclidean algorithm in the plane does not yet reveal the connection to the graph of the GCD. To that end, we turn the Euclidean algorithm on its head.

% by modifying our basis of the integer lattice $\ZZ^2$ in each step, instead of the point whose gcd we wish to compute.

\begin{figure}[t]
\begin{center}
\includegraphics[width=12cm]{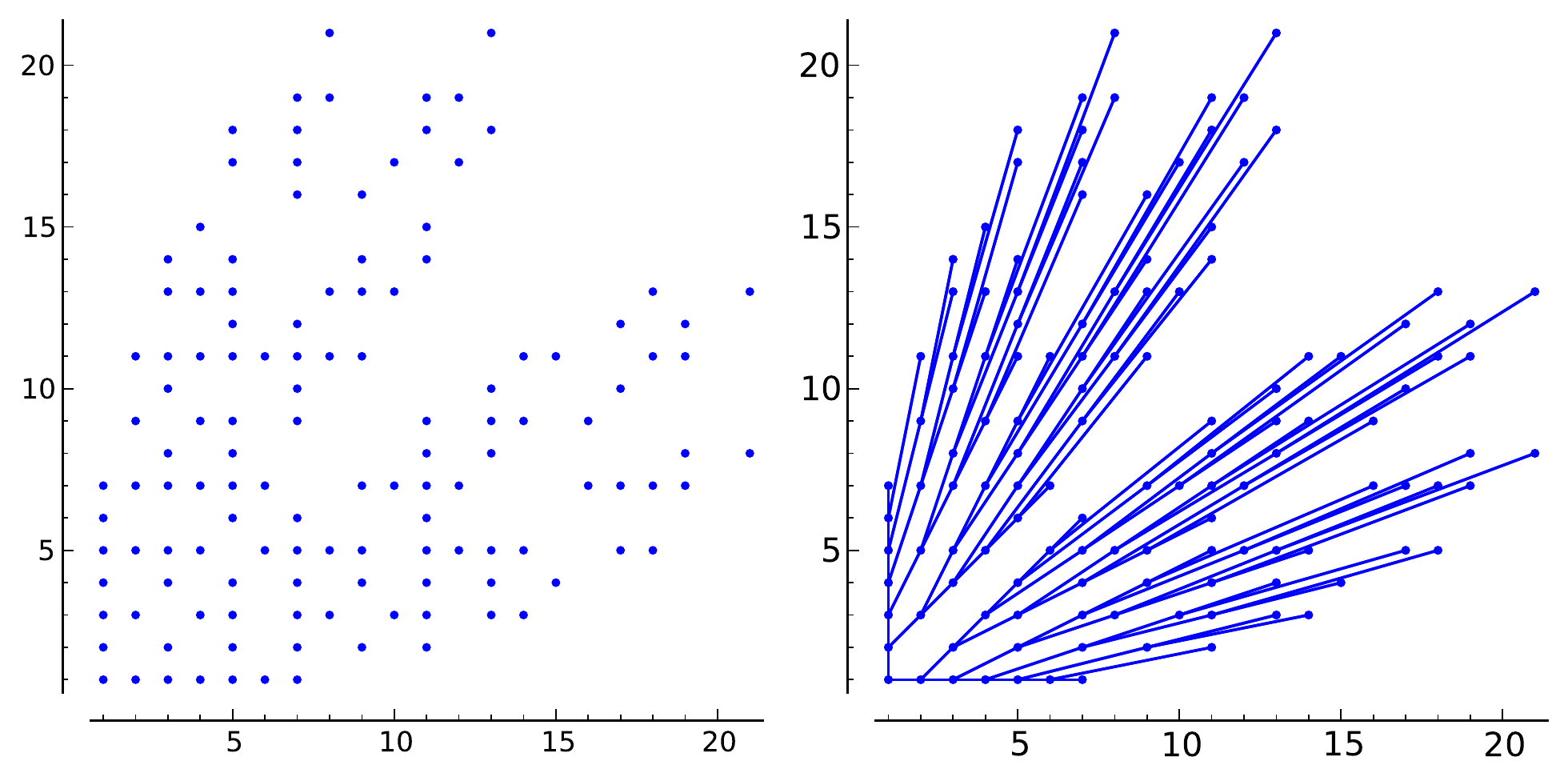}
\end{center}
\caption{The left panel shows all lattice points $(a,b)$ in the plane with $\gcd(a,b)=1$, up to a recursion depth of 5 in the Eulidean algorithm. The right panel also shows the tree structure induced by the inverted Euclidean algorithm as described in the text.}
\label{fig:euclid}
\end{figure}

We fix the point $p=(a,b)$ whose gcd we wish to compute and run the Euclidean algorithm by changing the basis $v_1$, $v_2$ of $\ZZ^2$ in each step. We define the \emph{center} of the current basis as the sum $c=v_1+v_2$. If $p$ lies below the line through $c$ we change our basis to $v_1'=v_1$ and $v_2'=c$. If $p$ lies above the line through $c$ we change our basis to $v_1'=c$ and $v_2'=v_2$. If $p$ lies on the line through $c$ we are done since $\gcd(a,b) = \frac{a}{c_1} = \frac{b}{c_2}$. Tracing out all the paths the center can take throughout this recursion, we obtain Figure~\ref{fig:euclid} which reveals the tree structure of the base points of the rays in Figure~\ref{fig:gcd-graph} and which gives a very natural (and novel) embedding of the Stern-Brocot tree \cite[p.~116-117]{Graham1994} in the plane.

\begin{figure}[t]
\begin{center}
\includegraphics[angle=270,width=12.2cm]{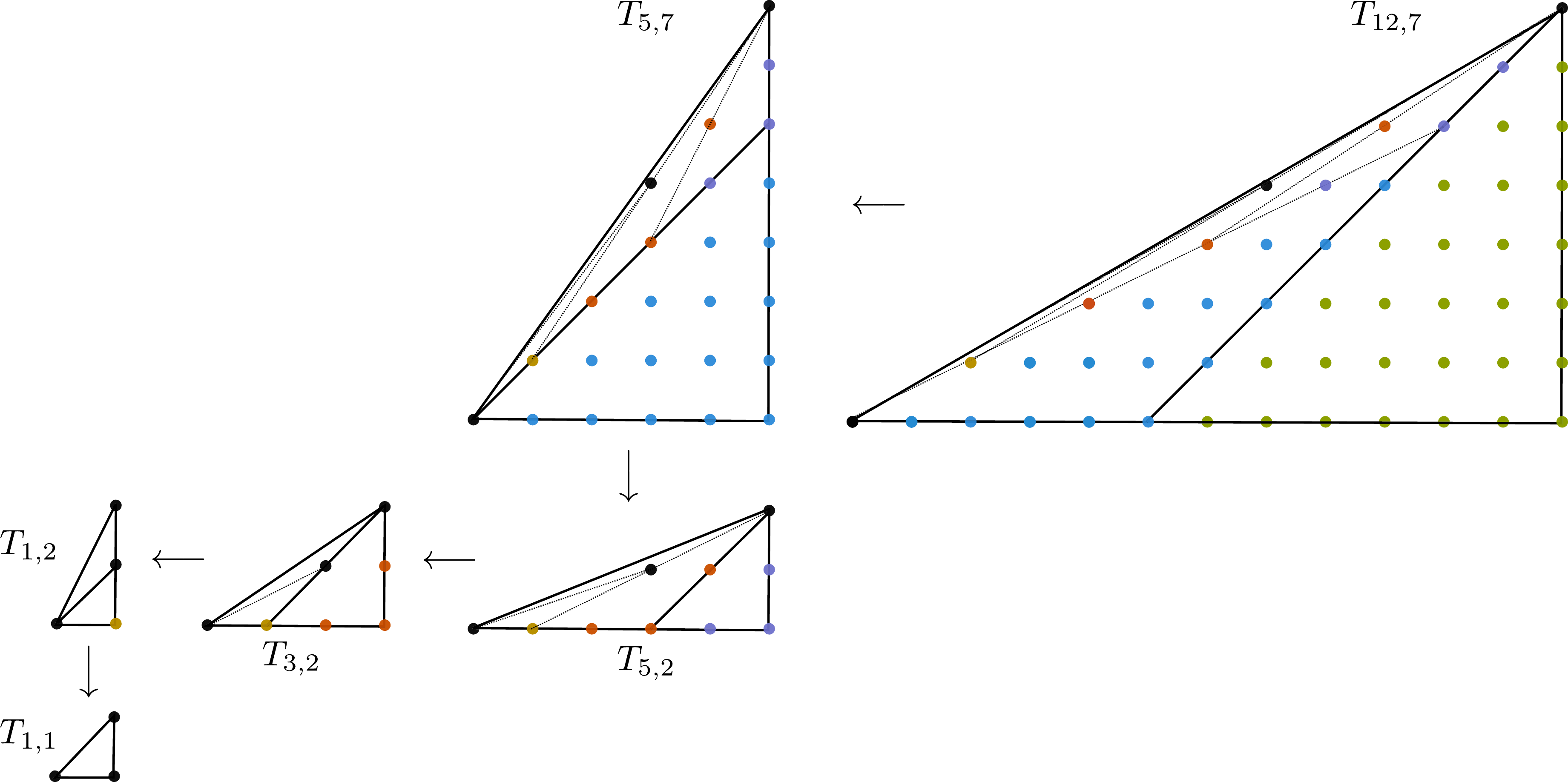}
\end{center}
\caption{Follwing the Euclidean algorithm, we reduce the triangle $T_{12,7}$ recursively, by removing triangles with integral slope and applying shearing lattice transformations. In this way we can decompose the ``staircase'' of integer points below the line from the origin to $(12,7)$ into ``simple'' triangles.}
\label{fig:staircase}
\end{figure}

To conclude, we follow \cite{BreuerHeymann2009} and examine the structure of the integer points below the line $L$ in more detail, going beyond the closest point $(p,q)$. Define $T_{a,b}$ to be the triangle with vertices $(0,0)$, $(a,0)$ and $(a,b)$. As we can see in Figure~\ref{fig:staircase}, the ``staircase'' of integer points in $T_{a,b}$ is irregular: the possible steps as we move from one column to the next are of two different heights, and it is not clear a priori what the underlying pattern is. It turns out, however, that the triangles $T_{a,b}$ have a very nice recursive structure. The key observation is that triangles of the form $T_{c,c}$ are very easy to describe as we always go exactly one step higher as we move from one column to the next. However, if $a > b$, then $T_{a,b}$ contains a triangle of the form $T_{b,b}$, sitting in the lower right corner. Removing this translate of the half-open triangle $T_{b,b}'$, we are left with a triangle $\tilde{T}$ with vertices $(0,0), (a-b,0)$ and $(a,b)$. Shearing $\tilde{T}$ using the linear transformation $A:(x,y)\mapsto (x-y,y)$, we see that the integer points in $\tilde{T}$ have the same structure as those in $T_{a-b,b}$. Here it is crucial that the linear transformation $A$ maps $\ZZ^2$ bijectively onto itself. Now, if $a<b$ we can apply the same procedure in the other direction. Just like in the Euclidean algorithm we continue recursively until we reach a triangle of the form $T_{c,c}$ at which point we stop. We can thus decompose any triangle $T_{a,b}$ into simple triangles of the form $T_{c,c}$. This process is illustrated in Figure~\ref{fig:staircase}.

This basic approach can yield much more information as detailed in \cite{BreuerHeymann2009}. For example, an analysis of how exactly the big and small steps in the staircase are distributed leads to several characterizations of Sturmian sequences of rational numbers. Moreover, using a recursive procedure similar in spirit to Figure~\ref{fig:staircase}, it is possible to show that the sets of lattice points in 2-dimensional fundamental parallelepipeds always have a short positive description as a union of Minkowski sums of discrete line segments -- this short description yields short rational function expressions for 2-dimensional fundamental parallelepipeds, which are quite distinct from those obtained via Barvinok's algorithm.

\paragraph{Acknowledgements.}
I would like to thank Benjamin Nill, Peter Paule, Manuel Kauers, Christoph Koutschan and an anonymous referee for their helpful comments on earlier versions of this article. I would also like to thank Matthias Beck whose lectures and book \cite{Beck2007} were my very own invitation to Ehrhart theory.

\bibliographystyle{splncs03}
\bibliography{references}

\begin{thebibliography}{10}
\providecommand{\url}[1]{\texttt{#1}}
\providecommand{\urlprefix}{URL }

\bibitem{Andrews2001}
Andrews, G., Paule, P., Riese, A.: {MacMahon's partition analysis VI: A new
  reduction algorithm}. Annals of Combinatorics  5(3),  251--270 (2001)

\bibitem{Monforte2013}
Aparicio~Monforte, A., Kauers, M.: {Formal Laurent Series in Several
  Variables}. Expositiones Mathematicae  31(4),  350--367 (2013)

\bibitem{Baldoni2011}
Baldoni, V., Berline, N., {De Loera}, J.A., K\"{o}ppe, M., Vergne, M.: {How to
  integrate a polynomial over a simplex}. Mathematics of Computation  80,
  297--325 (2011)

\bibitem{Barvinok2003}
Barvinok, A., Woods, K.: {Short rational generating functions for lattice point
  problems}. Journal of the American Mathematical Society  16(4),  957--979
  (2003)

\bibitem{Barvinok1994}
Barvinok, A.I.: {A Polynomial Time Algorithm for Counting Integral Points in
  Polyhedra When the Dimension Is Fixed}. Mathematics of Operations Research
  19(4),  769--779 (1994)

\bibitem{Barvinok2008}
Barvinok, A.I.: {Integer Points in Polyhedra}. European Mathematical Society
  (2008)

\bibitem{BBGM2014}
Beck, M., Breuer, F., Godkin, L., Martin, J.L.: {Enumerating colorings,
  tensions and flows in cell complexes}. Journal of Combinatorial Theory,
  Series A  122,  82--106 (Feb 2014)

\bibitem{Beck2005}
Beck, M., {De Loera}, J.A., Develin, M., Pfeifle, J., Stanley, R.P.:
  {Coefficients and roots of Ehrhart polynomials}. In: Barvinok, A.I. (ed.)
  Integer Points in Polyhedra, Proceedings of an AMS-IMS-SIAM Joint Summer
  Research Conference (Snowbird, Utah, 2003). pp. 1--24. AMS (2005)

\bibitem{Beck2009}
Beck, M., Haase, C., Sottile, F.: {(Formulas of Brion, Lawrence, and Varchenko
  on rational generating functions for cones)}. The Mathematical Intelligencer
  31(1),  9--17 (2009)

\bibitem{Beck2007}
Beck, M., Robins, S.: {Computing the continuous discretely: Integer-point
  enumeration in polyhedra}. Springer (2007)

\bibitem{Beck2014}
Beck, M., Sanyal, R.: {Combinatorial reciprocity theorems} (2014),
  \url{http://math.sfsu.edu/beck/crt.html}, to appear

\bibitem{Beck2006-iop}
Beck, M., Zaslavsky, T.: {Inside-out polytopes}. Advances in Mathematics
  205(1),  134--162 (2006)

\bibitem{Beck2006-flow}
Beck, M., Zaslavsky, T.: {The number of nowhere-zero flows on graphs and signed
  graphs}. Journal of Combinatorial Theory, Series B  96(6),  901--918 (2006)

\bibitem{Bergeron2009}
Bergeron, N., Zabrocki, M.: {The Hopf algebras of symmetric functions and
  quasi-symmetric functions in non-commutative variables are free and co-free}.
  Journal of Algebra and Its Applications  8(4),  581--600 (2009)

\bibitem{Breuer2009}
Breuer, F.: {Ham Sandwiches, Staircases and Counting Polynomials}. Phd thesis,
  Freie Universit\"{a}t Berlin (2009)

\bibitem{Breuer2012}
Breuer, F.: {Ehrhart $f^*$-coefficients of polytopal complexes are non-negative
  integers}. Electronic Journal of Combinatorics  19(4),  P16 (2012)

\bibitem{Breuer2011}
Breuer, F., Dall, A.: {Bounds on the Coefficients of Tension and Flow
  Polynomials}. Journal of Algebraic Combinatorics  33(3),  465--482 (2011)

\bibitem{BDK2012}
Breuer, F., Dall, A., Kubitzke, M.: {Hypergraph coloring complexes}. Discrete
  mathematics  312(16),  2407--2420 (Aug 2012)

\bibitem{BEKZ14}
Breuer, F., Eichhorn, D., Kronholm, B.: Cranks and the geometry of
  combinatorial witnesses for the divisibility and periodicity of the
  restricted partition function  (2014), in preparation

\bibitem{BreuerHeymann2009}
Breuer, F., von Heymann, F.: {Staircases in $\mathbb{Z}^2$}. Integers  10(6),
  807--847 (2010)

\bibitem{Breuer2014}
Breuer, F., Klivans, C.J.: {Scheduling problems}  (2014), submitted,
  arXiv:1401.2978v1

\bibitem{BreuerSanyal2012}
Breuer, F., Sanyal, R.: {Ehrhart theory, Modular flow reciprocity, and the
  Tutte polynomial}. Mathematische Zeitschrift  270(1),  1--18 (2012)

\bibitem{BZ14}
Breuer, F., Zafeirakopoulos, Z.: {Polyhedral Omega}: A new algorithm for
  solving linear {Diophantine} systems  (2014), in preparation

\bibitem{Brion1988}
Brion, M.: Points entiers dans les poly{\`e}dres convexes. Annales
  scientifiques de l'{\'E}cole Normale Sup{\'e}rieure  21(4),  653--663 (1988)

\bibitem{BIS2012}
Bruns, W., Ichim, B., S\"oger, C.: The power of pyramid decomposition in
  {Normaliz}  (2012), arXiv:1206.1916v1

\bibitem{Calkin2000}
Calkin, N., Wilf, H.S.: {Recounting the rationals}. The American Mathematical
  Monthly  107(4),  360--363 (2000)

\bibitem{Chari1997}
Chari, M.K.: {Two decompositions in topological combinatorics with applications
  to matroid complexes}. Transactions of the American Mathematical Society
  349(10),  3925--3943 (1997)

\bibitem{DeLoera2009}
{De Loera}, J.A., Hemmecke, R., K\"{o}ppe, M.: {Pareto Optima of Multicriteria
  Integer Linear Programs}. INFORMS Journal on Computing  21(1),  39--48 (2009)

\bibitem{DeLoera2012}
De~Loera, J.A., Hemmecke, R., K{\"o}ppe, M.: Algebraic and geometric ideas in
  the theory of discrete optimization, vol.~14. SIAM (2012)

\bibitem{DeLoera2004}
{De Loera}, J.A., Hemmecke, R., Tauzer, J., Yoshida, R.: {Effective lattice
  point counting in rational convex polytopes}. Journal of Symbolic Computation
   38(4),  1273--1302 (Oct 2004)

\bibitem{Ehrhart1962}
Ehrhart, E.: {Sur les poly\`{e}dres rationnels homoth\'{e}tiques \`{a} $n$
  dimensions}. C. R. Acad. Sci. Paris  254,  616--618 (1962)

\bibitem{Fukuda1996}
Fukuda, K., Prodon, A.: Double description method revisited. In: Deza, M.,
  Euler, R., Manoussakis, I. (eds.) Combinatorics and Computer Science, Lecture
  Notes in Computer Science, vol. 1120, pp. 91--111. Springer Berlin Heidelberg
  (1996)

\bibitem{Fukuda1994}
Fukuda, K., Rosta, V.: {Combinatorial face enumeration in convex polytopes}.
  Computational Geometry  4(4),  191--198 (Aug 1994)

\bibitem{Graham1994}
Graham, R.L., Knuth, D.E., Patashnik, O.: Concrete Mathematics: A Foundation
  for Computer Science. Addison-Wesley Longman Publishing Co., Inc., Boston,
  MA, USA, 2nd edn. (1994)

\bibitem{Greene1977}
Greene, C.: {Acyclic orientations (Notes)}. Higher Combinatorics (M. Aigner,
  ed.), Reidel, Dordrecht pp. 65--68 (1977)

\bibitem{Haase2009}
Haase, C., Schicho, J.: {Lattice polygons and the number 2i+7}. The American
  Mathematical Monthly  116(2),  151--165 (2009)

\bibitem{Henk2009}
Henk, M., Tagami, M.: {Lower bounds on the coefficients of Ehrhart
  polynomials}. European Journal of Combinatorics  30(1),  70--83 (Jan 2009)

\bibitem{Hersh2007}
Hersh, P., Swartz, E.: {Coloring complexes and arrangements}. Journal of
  Algebraic Combinatorics  27(2),  205--214 (Jul 2007)

\bibitem{Jochemko2013}
Jochemko, K., Sanyal, R.: {Arithmetic of marked order polytopes, monotone
  triangle reciprocity, and partial colorings} pp. 1--16 (2013),
  arXiv:1206.4066v2

\bibitem{Koppe2008}
K\"{o}ppe, M., Verdoolaege, S.: {Computing parametric rational generating
  functions with a primal Barvinok algorithm}. Electronic Journal of
  Combinatorics  15,  R16 (2008)

\bibitem{Lawrence1988}
Lawrence, J.: Valuations and polarity. Discrete \& Computational Geometry
  3(1),  307--324 (1988)

\bibitem{Lenstra1983}
Lenstra, H.W.: {Integer Programming with a Fixed Number of Variables}.
  Mathematics of Operations Research  8(4),  538--548 (1983)

\bibitem{Macdonald1971}
Macdonald, I.G.: {Polynomials associated to finite cell complexes}. Journal of
  the London Mathematical Society (2)  4,  181--192 (1971)

\bibitem{Pfeifle2003}
{Pfeifle}, J., {Rambau}, J.: {Computing triangulations using oriented
  matroids.} In: {Algebra, geometry, and software systems}, pp. 49--75. Berlin:
  Springer (2003)

\bibitem{DeLoera2010}
{Rambau}, J.A.D.J., {Santos}, F.: Triangulations. Structures for algorithms and
  applications. Berlin: Springer (2010)

\bibitem{Schrijver1986}
Schrijver, A.: {Theory of Linear and Integer Programming}. John Wiley
  $\backslash$\& Sons, Inc. (1986)

\bibitem{Stanley1973}
Stanley, R.P.: {Acyclic orientations of graphs}. Discrete Mathematics  5,
  171--178 (May 1973)

\bibitem{Stanley1980}
Stanley, R.P.: {Decompositions of Rational Convex Polytopes}. Annals of
  Discrete Mathematics  6,  333--342 (1980)

\bibitem{Stanley1986}
Stanley, R.P.: {Two poset polytopes}. Discrete \& Computational Geometry  1,
  9--23 (1986)

\bibitem{Stanley2001}
Stanley, R.P.: Enumerative Combinatorics, vol.~2. Cambridge Studies in Advanced
  Mathematics (2001)

\bibitem{Stapledon2009}
Stapledon, A.: Inequalities and {Ehrhart} $\delta$-vectors. Transactions of the
  American Mathematical Society  361,  5615--5626 (2009)

\bibitem{Swartz2006}
Swartz, E.: {g-Elements, finite buildings and higher Cohen--Macaulay
  connectivity}. Journal of Combinatorial Theory, Series A  113,  1305--1320
  (2006)

\bibitem{Varchenko1987}
Varchenko, A.N.: Combinatorics and topology of the disposition of affine
  hyperplanes in real space. Functional Analysis and Its Applications  21(1),
  9--19 (1987)

\bibitem{Verdoolaege2008}
Verdoolaege, S., Woods, K.: {Counting with rational generating functions}.
  Journal of Symbolic Computation  43(2),  75--91 (2008)

\bibitem{Verdoolaege2007}
Verdoolaege, S., Seghir, R., Beyls, K., Loechner, V., Bruynooghe, M.: {Counting
  Integer Points in Parametric Polytopes Using Barvinok's Rational Functions}.
  Algorithmica  48(1),  37--66 (2007)

\bibitem{Woods2013}
Woods, K.: {Presburger Arithmetic, Rational Generating Functions, and
  Quasi-Polynomials}. In: Fomin, F.V., Freivalds, R., Kwiatkowska, M., Peleg,
  D. (eds.) Automata, Languages, and Programming. pp. 410--421. Springer Berlin
  Heidelberg (2013), \url{http://arxiv.org/abs/1211.0020}

\bibitem{Ziegler1995}
Ziegler, G.M.: {Lectures on Polytopes}. Graduate Texts in Mathematics, Springer
  (1995)

\end{thebibliography}

\end{document}